\documentclass[letterpaper]{article}

\date{\today}

\newcommand{\myauthor}{Benjamin Antieau}
\newcommand{\mytitle}{\v{C}ech approximation to the Brown-Gersten spectral sequence}
\newcommand{\pdftitle}{Cech approximation to the Brown-Gersten spectral sequence}

\title{\mytitle}
\author{\myauthor\footnote{This material is based upon work supported by the NSF under Grant No. DMS-0901373.}}

\usepackage[pdfstartview=FitH,
            pdfauthor={\myauthor},
            pdftitle={\pdftitle},
            colorlinks,
            linkcolor=reference,
            citecolor=citation,
            urlcolor=e-mail,
            backref]{hyperref}

\usepackage{rnote}
\usepackage{theorems}
\usepackage{times}
\usepackage[all,cmtip]{xy}
\usepackage{diagxy}
\usepackage{tikz}
\usetikzlibrary{matrix,arrows}

\begin{document}
\maketitle

\begin{abstract}
  \noindent
    In this paper, we show that the \'etale index of a torsion cohomological Brauer class is divisible by the period of the class. The tool used to make this computation is
    the \v{C}ech approximation of the title. To create the approximation, we use the folklore theorem that the homotopy limit and Postnikov spectral sequences for a cosimplicial space
    agree beginning with the $\Eoh_2$-page. As far we know, this folklore theorem has no proof in the literature, so we include a proof.


\paragraph{Mathematics Subject Classification 2010}
Primary: \href{http://www.ams.org/mathscinet/msc/msc2010.html?t=14Fxx&btn=Current}{14F22}, \href{http://www.ams.org/mathscinet/msc/msc2010.html?t=18Gxx&btn=Current}{18G40}.
Secondary: \href{http://www.ams.org/mathscinet/msc/msc2010.html?t=16Kxx&btn=Current}{16K50}, \href{http://www.ams.org/mathscinet/msc/msc2010.html?t=19Dxx&btn=Current}{19D23}.
\end{abstract}

\tableofcontents

\section{Introduction}
Let $X$ be a geometrically connected scheme, and let $\Br'(X)=\Hoh^2(X_{\et},\Gm)_{\tors}$ be the cohomological Brauer group of $X$.
There are two integer invariants of a Brauer class $\alpha\in\Br'(X)$. The first is the period of $\alpha$;
it is the order of $\alpha$ in the group $\Br'(X)$, and it is generally written as $per(\alpha)$. The second is the index of $\alpha$. If $\mathcal{A}$ is an Azumaya algebra on $X$, then
it has rank $n^2$ as a locally free $\mathcal{O}_X$ sheaf. The integer $n$ is called the index of $\mathcal{A}$. For a class $\alpha$, its index $ind(\alpha)$ is defined as the least integer $n$
such that there is an Azumaya algebra $\mathcal{A}$ in the class of $\alpha$ of index $n$. If no such Azumaya algebra exists, then the index is defined to be $+\infty$.

In general,
\begin{equation}\label{eq:div1}
    per(\alpha)|ind(\alpha).
\end{equation}
If $X=\Spec k$ for a field $k$, then the index is finite, and the period and the index have the same prime divisors. Therefore, there is a least integer $e(\alpha)$ such that
\begin{equation*}
    ind(\alpha)|(per(\alpha))^{e(\alpha)}.
\end{equation*}
As $\alpha$ ranges over all Brauer classes, it is an interesting open question to determine the values that $e(\alpha)$ can take on. In \cite{colliot_brauersche_2001}, \colliot asks whether
the following is true.

\begin{conjecture}[Period-Index Conjecture]\label{conj:perind}
    Let $k$ be a field of dimension $d>0$. Then,
    \begin{equation}\label{eq:perindconj}
        ind(\alpha)|(per(\alpha))^{d-1},
    \end{equation}
    for all $\alpha\in\Br(k)$.
\end{conjecture}

The notion of dimension in the conjecture is not entirely settled. The conjecture is made specifically~\cite{colliot_brauersche_2001} for function fields of $d$-dimensional varieties
over algebraically closed fields, or for $C_d$ fields. But it is known to be false for cohomological dimension.

Several cases in low dimensions are known, but the question is open in general. For example, it is not known for any function field $\CC(X)$ of an algebraic $3$-fold $X$.
The importance of the conjecture for $k(X)$ is that it gives information both about the geometry of $X$ and about the arithmetic of $k(X)$ when $X$ is an algebraic variety over an algebraically closed field $k$.

In \cite{antieau_index_2009}, we introduced a new invariant $eti(\alpha)$ for $\alpha\in\Br'(X)$, which we call the \'etale index.
By the definition,
\begin{equation*}
    eti(\alpha)|ind(\alpha).
\end{equation*}
The main theorem of that paper is an analogue of Equation \eqref{eq:perindconj} when $X=\Spec k$. Specifically, if $k$ is a field of finite cohomological dimension $d=2c$ or $d=2c+1$, then
\begin{equation}\label{eq:pereti}
    eti(\alpha)|(per(\alpha))^c
\end{equation}
when the prime divisors of $per(\alpha)$ are ``large'' with respect to $d$. Specifically, the statement is true when $d<2q$ for all primes $q$ that divide $per(\alpha)$.

The main result of this paper, Theorem \ref{thm:main}, is that the analogue of Equation \eqref{eq:div1} holds for the \'etale index as well:
\begin{equation*}
    per(\alpha)|eti(\alpha).
\end{equation*}
In particular, it is not equal to $1$ if $\alpha$ is non-trivial.
Therefore, the \'etale index possesses similar formal properties to the index. However, the \'etale index is always finite for schemes $X$ of finite \'{e}tale cohomological dimension, while
certain Brauer classes on non-quasi-projective schemes have infinite index. 

Here is an application of Theorem~\ref{thm:main}. Define
\begin{align*}
    \K_0^{\alpha}(k)^{(0)}  &=  \K_0^{\alpha}(k)/\ker\left(\rank:\K_0^{\alpha}(k)\rightarrow\ZZ\right)\\
    \K_0^{\alpha,\et}(k)^{(0)}  &=  \K_0^{\alpha,\et}(k)/\ker\left(\rank:\K_0^{\alpha,\et}(k)\rightarrow\ZZ\right),
\end{align*}
where $\K^{\alpha}$ is twisted connective $K$-theory, and $\K^{\alpha,\et}$ is the \'etale sheafification.
The image of the rank on $\K_0^{\alpha}(k)$ is generated by $ind(\alpha)$, while the image of the rank on $\K_0^{\alpha,\et}$ is generated by $eti(\alpha)$. There is an inclusion
\begin{equation}\label{eq:kinclusion}
    \K_0^{\alpha}(k)^{(0)}\subseteq\K_0^{\alpha,\et}(k)^{(0)}
\end{equation}
with cokernel $F=\ZZ/\left(\frac{ind(\alpha)}{eti(\alpha)}\right)$. In the untwisted case ($\alpha=0$), the cokernel is always zero. On the other hand, Equation~\eqref{eq:pereti}, together with the sharpness of some known cases of the
period-index conjecture, imply that $F$ is in general not zero when $\alpha\neq 0$.

Suppose that the period-index conjecture of \colliot is true. Then, the fact that $per(\alpha)|eti(\alpha)$ implies that the cokernel $F$ is of order \emph{at most} $per(\alpha)^{d-2}$.
Conversely, if one could prove that the cokernel was of order at most $per(\alpha)^{c}$ when $d$ is odd or $per(\alpha)^{c-1}$ when $d$ is even,
then Equation~\eqref{eq:pereti} would imply the period-index conjecture, at least away from some small primes.
Thus, our result allows a translation of the period-index conjecture into a question about the arithmetic failure of \'etale descent for twisted $K$-theory.

The proof of Theorem~\ref{thm:main} is based on the computation of a differential in the descent spectral sequence associated to $\K^{\alpha}$. Recall that the descent spectral sequence,
or Brown-Gersten spectral sequence, is
\begin{equation*}
    \Eoh_2^{s,t}=\Hoh^s(X_{\et},\mathcal{K}_t^{\alpha})\Rightarrow\K_{t-s}^{\alpha,\et}(X),
\end{equation*}
with differentials $d_r^{\alpha}$ of degree $(r,r-1)$.
In \cite{antieau_index_2009}, we proved that the twisted $K$-theory sheaves are isomorphic, in a natural way, to the untwisted $K$-theory sheaves. In particular, $\mathcal{K}_0^{\alpha}\iso\ZZ$,
and $\mathcal{K}_1^{\alpha}\iso\Gm$. The \'etale index is defined to be the least integer $n$ such that $d_r^{\alpha}(n)=0$, $r\geq 2$, where we view $n$ as an element in $\Hoh^0(X,\mathcal{K}_0^{\alpha})$.
The proof of the main theorem, that $per(\alpha)|eti(\alpha)$, is given by proving that $d_2^{\alpha}(1)=\alpha$, where $\alpha\in\Hoh^2(X_{\et},\Gm)\iso\Hoh^2(X_{\et},\mathcal{K}_1^{\alpha})$.
Indeed, once this is done, the least integer such that $d_2^{\alpha}(n)=0$ is $n=per(\alpha)$. This computation is analogous to
one given by Kahn and Levine in \cite[Proposition~6.9.1]{kahn_levine} for the \'etale motivic spectral sequence converging to \'etale twisted $K$-theory
and another given by Atiyah and Segal in \cite[Proposition~4.6]{atiyah_twisted_2006} in the topological twisted case.

To make the computation of the differential, we introduce a \v{C}ech approximation to the descent spectral sequence. This in turn relies on the comparison theorem, Theorem~\ref{thm:isomorphism},
which says that the homotopy limit and Postnikov spectral sequences for cosimplicial spaces agree beginning with the $\Eoh_2$-page. This is a folklore theorem, and we include a proof since we know
of no reference for it. Both types of spectral sequences arise frequently in practice in homotopy theory due to the mechanism of cosimplicial replacement of presheaves of spaces, and the theorem easily
applies in this situation to say that for a presheaf of spaces on a small category, the homotopy limit and Postnikov spectral sequences agree beginning with the $\Eoh_2$-page.

Once this approximation is established, it remains to describe the $d_2$ differential for the homotopy limit spectral sequence of a cosimplicial space, and then to translate
this description to the \v{C}ech approximation spectral sequence.

Some background is given in Sections~\ref{sub:towers} and \ref{sec:cosimplicial}. Several spectral sequences for cosimplicial spaces are introduced, and a comparison theorem proven, in Section~\ref{sec:ss}.
This is applied to presheaves of spectra on a Grothendieck
site in Section~\ref{sec:cech}. The differentials of the homotopy limit spectral sequence are described in Section~\ref{sec:differentials},
and those of the \v{C}ech approximation in Section~\ref{sec:cechd2}.
The main theorem is proven in Section~\ref{sec:division}. Finally, in Section~\ref{sec:indeti}, the cokernel $F$ of Equation~\eqref{eq:kinclusion} is related in detail to the period-index conjecture.

\paragraph{Acknowledgments}
This paper is part of my Ph.D. dissertation, which was written under the direction of Henri Gillet at UIC. I owe him a great deal of thanks for his support.
I have also spoken with Peter Bousfield, David Gepner, Christian Haesemeyer, Jumpei Nogami, and Brooke Shipley about this work. They have all been very helpful and encouraging.
Additionally, the referee did a very careful job which certainly has resulted in an improvement in the exposition.

\section{Spectral Sequences Associated to Towers of Fibrations}\label{sub:towers}

The first quadrant spectral sequences of this section are generalized in the sense that not every term is an abelian group.
They have the property $\Eoh_1^{s,t}$ is a group when $t-s=1$ and it is a pointed set when $t-s=0$.
Nonetheless, there is a good notion of such spectral sequences. For details on these generalized spectral sequences as well as all the other material and notation in this section, see
\cite[Section~IX.4, esp. 4.1]{bousfield_homotopy_1972}.

\begin{construction}
    Let
    \begin{equation*}
        \cdots\rightarrow X_n\rightarrow X_{n-1}\rightarrow\cdots\rightarrow X_1\rightarrow X_0\rightarrow\ast
    \end{equation*}
    be a tower of fibrations of pointed spaces. Let $F_s$ be the fiber of $X_s\rightarrow X_{s-1}$, and denote by $X$ the limit of the tower.

    There are long exact sequences associated to the fibrations $F_s\rightarrow X_s\rightarrow X_{s-1}$:
    \begin{equation}\label{eq:lowerterms}
        \cdots\rightarrow\pi_{t-s}F_s\xrightarrow{k}\pi_{t-s}X_s\xrightarrow{i}\pi_{t-s}X_{s-1}\xrightarrow{j}\pi_{t-s-1}F_s\rightarrow\cdots.
    \end{equation}
    These sequences continue all the way down to $\pi_0 X_{s-1}$:
    \begin{equation*}
        \pi_2 X_{s-1}\xrightarrow{j}\pi_1 F_s\xrightarrow{k}\pi_1 X_s\xrightarrow{i}\pi_1 X_{s-1}\xrightarrow{j}\pi_0 F_s\xrightarrow{k}\pi_0 X_s\xrightarrow{i}\pi_0 X_{s-1},
    \end{equation*}
    and
    \begin{equation*}
        \pi_1 X_{s-1}\xrightarrow{j}\pi_0 F_s
    \end{equation*}
    extends to an action of $\pi_1 X_{s-1}$ on $\pi_0 F_s$ so that $j$ is the map onto the orbit of the basepoint of $F_s$ under the action of $\pi_1 X_{s-1}$. Besides the usual
    conditions of $\ker=\im$ in the range that this makes sense, the exactness
    of Equation~\eqref{eq:lowerterms} \emph{means}
    \begin{itemize}
        \item   that the quotient of $\pi_0 F_s$ under this action injects into $\pi_0 X_s$,
        \item   that the cokernel (quotient set) of $\pi_0 F_s\xrightarrow{k}\pi_0 X_s$ injects into $\pi_0 X_{s-1}$,
        \item   that the stabilizer of the action of $\pi_1 X_{s-1}$ on $\pi_0 F_s$ at the base-point of $F_s$ is the quotient of $\pi_1 X_s$ by the image of $\pi_1 F_s\rightarrow \pi_1 X_s$, and
        \item   that $\pi_2 X_{s-1}$ maps to the center of $\pi_1 F_s$.
    \end{itemize}

    The tower of fibrations and the exact sequences above define an exact couple
    \begin{center}
        \begin{tikzpicture}[description/.style={fill=white,inner sep=2pt}]
                \matrix (m) [matrix of math nodes, row sep=3em,
                column sep=2.5em, text height=1.5ex, text depth=0.25ex]
                {     D_1 & & D_1\\
                      & E_1 &\\};
                \path[->,font=\scriptsize]
                (m-1-1) edge node[auto] {$ i $} (m-1-3)
                (m-1-3) edge node[below right] {$ j $} (m-2-2)
                (m-2-2) edge node[below left] {$ k $} (m-1-1);
        \end{tikzpicture}
    \end{center}
    where $D$ and $E$ are bigraded groups:
    \begin{align*}
        \Doh_1^{s,t}   &=  \pi_{t-s}X_s\\
        \Eoh_1^{s,t}   &=  \pi_{t-s}F_s.
    \end{align*}
    The maps $i$, $j$, and $k$ are of bi-degrees $(-1,1)$, $(1,-2)$, and $(0,0)$:
    \begin{align*}
        i:  \pi_{t-s}X_s      &\rightarrow    \pi_{t-s}X_{s-1}\\
        j:  \pi_{t-s}X_{s-1}  &\rightarrow    \pi_{t-s-1}F_s\\
        k:  \pi_{t-s-1}F_s    &\rightarrow    \pi_{t-s-1}X_s.
    \end{align*}
    
    As usual, the exact couple gives rise to a differential $d=j\circ k$ on $E$. It is of bi-degree $(1,-2)$.
    The first derived exact couple is
    \begin{align*}
        \pi_{t-s}X_s^{(1)}:=\Doh_2^{st} &=  \im(i)=\im(\pi_{t-s}X_{s+1}\xrightarrow{i}\pi_{t-s}X_s)\subseteq\pi_{t-s}X_s\\
        \pi_{t-s}F_s^{(1)}:=\Eoh_2^{st} &=  \Hoh(d)=\ker(\pi_{t-s}F_s\xrightarrow{k}\pi_{t-s}X_s/\im(i))/\\
                                       &   \ker(\pi_{t-s+1}X_{s-1}\xrightarrow{i}\pi_{t-s+1}X_{s-2}).
    \end{align*}
    When $s=t$, then the definition of $\Eoh_2^{st}$ should be interpreted as the quotient of the pointed set
    \begin{equation*}
        \ker(\pi_0 F_s\xrightarrow{k}\pi_0 X_s/\im(i))
    \end{equation*}
    by the action of $\ker(\pi_1 X_{s-1}\xrightarrow{i}\pi_1 X_{s-2})\subseteq\pi_1 X_{s-1}$. Then, the sequences
    \begin{equation*}
    \pi_2 X_{s-1}^{(1)}\xrightarrow{j}\pi_1 F_s^{(1)}\xrightarrow{k}\pi_1 X_s^{(1)}\xrightarrow{i}\pi_1 X_{s-1}^{(1)}\xrightarrow{j}\pi_0 F_s^{(1)}\xrightarrow{k}\pi_0 X_s^{(1)}\xrightarrow{i}\pi_0 X_{s-1}^{(1)},
    \end{equation*}
    are also exact in the generalized sense above.

    Repeating this process, one obtains a generalized spectral sequence $\Eoh_1^{s,t}\{X_\ast\}$ associated to the tower, with
    \begin{equation*}
        \Eoh_1^{s,t}=\pi_{t-s} F_s\rightharpoonup\pi_{t-s} X,
    \end{equation*}
    where $F_s$ is the fiber of $X_n\rightarrow X_{n-1}$. The differential $d_r$ is of degree $(r,r-1)$. See \cite[Chapter~IX]{bousfield_homotopy_1972}.
    The harpoon $\rightharpoonup$ means that the spectral sequence may not converge in the usual sense. Instead, there is a filtration
    \begin{equation*}
        Q_s\pi_i X=\ker(\pi_i X\rightarrow\pi_i X_s)
    \end{equation*}
    with successive quotients
    \begin{equation*}
        e_{\infty}^{s,t}=\ker(Q_s\pi_{t-s} X\rightarrow Q_{s-1}\pi_{t-s} X),
    \end{equation*}
    and inclusions
    \begin{equation}\label{eq:convergence}
        e_{\infty}^{s,t}\subseteq\Eoh_{\infty}^{s,t}.
    \end{equation}
    Write $\Eoh_1^{s,t}\Rightarrow\pi_{t-s} X$ when the spectral sequence does hold in the usual sense, in which case equality holds in Equation~\eqref{eq:convergence}. In this case,
    say that the spectral sequence \df{converges completely}.

    There is a second spectral sequence $\tilde{\Eoh}_{2}^{s,t}\{X_\ast\}$, which is simply a re-indexed version of the first:
    \begin{equation*}
        \tilde{\Eoh}_2^{s,t}=\Eoh_1^{t,2t-s}=\pi_{t-s} F_{t}\rightharpoonup \pi_{t-s} X.
    \end{equation*}
    This is derived from the exact couple
    \begin{center}
        \begin{tikzpicture}[description/.style={fill=white,inner sep=2pt}]
                \matrix (m) [matrix of math nodes, row sep=3em,
                column sep=2.5em, text height=1.5ex, text depth=0.25ex]
                {     \tilde{D}_2 & & \tilde{D}_2\\
                & \tilde{E}_2, &\\};
                \path[->,font=\scriptsize]
                (m-1-1) edge node[auto] {$ i $} (m-1-3)
                (m-1-3) edge node[below right] {$ j $} (m-2-2)
                (m-2-2) edge node[below left] {$ k $} (m-1-1);
        \end{tikzpicture}
    \end{center}
    With $\tilde{\Eoh}_2$ as above, and
    \begin{equation*}
        \tilde{\Doh}_2^{s,t}=\pi_{t-s}X_t.
    \end{equation*}
    The filtration is the same:
    \begin{equation*}
        \tilde{Q}_s\pi_i X=Q_s\pi_i X,
    \end{equation*}
    but the successive quotients are
    \begin{equation*}
        \tilde{e}_{\infty}^{s,t}=e_{\infty}^{t,2t-s}.
    \end{equation*}
    Here, the differentials are also of degree $(r,r-1)$.

    The spectral sequence and the filtration $Q_s\pi_i$ are functorial for towers of pointed spaces. We view the spectral sequence as including the information of the abutment and
    the filtration on the abutment. Thus, a morphism of spectral sequences includes a filtration-respecting morphism of the abutment.
\end{construction}

\begin{remark}
    Under certain conditions, these spectral sequence do converge in some range to some of the homotopy groups of $X$.
    For instance, suppose that $i\geq1$ and that for each $s\geq 0$ there is an integer $N(s)\geq 1$ such that
    \begin{equation}\label{eq:convergence}
        \Eoh_{M}^{s,s+j}=\Eoh_{\infty}^{s,s+j}
    \end{equation}
    for all $M\geq N(s)$ when $j=i$ and $j=i+1$. Then, $\Eoh_{\ast}^{s,t}\{X_\ast\}$ converges to $\pi_i X$. This is also true for $i=0$ when all of the homotopy sets of the spaces in the tower
    are abelian groups. Again, for details see \cite[Section~IX.5]{bousfield_homotopy_1972}. 
\end{remark}

\begin{remark}
    Note that in the application to $K$-theory, the spaces in the tower will have homotopy sets $\pi_t$ which are abelian groups for all $t\geq 0$.
    Moreover, the convergence conditions of Equation~\eqref{eq:convergence} will always hold under the finite cohomological dimension conditions used in this paper.
\end{remark}

\section{Cosimplicial Spaces}\label{sec:cosimplicial}

\begin{definition}
  Let $\Delta$ be the category of finite simplices. Objects of $\Delta$ are non-empty finite ordered sets, and morphisms are set morphisms that preserve order.
  The category $\sSets$ of simplicial sets is the functor category $\Fun(\Delta^{\op},\Sets)$. In general, if $C$ is a category, then $\mathbf{s}C$ is the category $\Fun(\Delta^{\op},C)$,
  the category of simplicial objects in $C$. Objects of $\sSets$ will be called spaces. The category $\sSets_\ast$ is the category of pointed spaces.
\end{definition}

\begin{definition}
  If $C$ is a category, then denote by $\mathbf{c}C$ the functor category $\Fun(\Delta,C)$, the category of cosimplicial objects in $C$. The category of cosimplicial spaces
  is the category $\csSets$. Write $\csSets_\ast$ for the category of cosimplicial pointed spaces.
\end{definition}

\begin{example}
    For a space $X$, let the same symbol $X$ denote the constant cosimplicial space $n\mapsto X$.
\end{example}

\begin{example}
    The cosimplicial space $\Delta$ is the functor $n\mapsto\Delta^n$, where $\Delta^n$ is the simplicial space $\Delta^n$.
\end{example}

\begin{example}
    Let $\mathcal{U}^{\bullet}$ be a hypercover in a site $C$, and let $X$ be a presheaf of simplicial sets on $C$ (a presheaf of spaces). Then, $X_{\mathcal{U}^{\bullet}}$ denotes the cosimplicial space given
    by evaluating $X$ at each level of $\mathcal{U}^{\bullet}$ in the usual way.
\end{example}

\begin{definition}
  If $F$ is an endofunctor of $\sSets$, then one extends $F$ to an endofunctor on $\csSets$ by level-wise application. That is, for a cosimplical space $X$, define $F(X)^n=F(X^n)$. 
  The typical examples are the $s$-skeleton functors $X\mapsto X[s]$ and the $\Ex$-functor.
\end{definition}

\begin{definition}
    Let $\mathbf{P}$ be a property of spaces. Then, a cosimplicial space $X$ is level $\mathbf{P}$ if each space $X_n$ is $\mathbf{P}$ for $n\geq 0$.
    Similarly, if $\mathbf{Q}$ is a property of morphisms of spaces, then a morphism $f:X\rightarrow Y$ of cosimplicial spaces is level $\mathbf{Q}$ if $f^n:X^n\rightarrow Y^n$ is $\mathbf{Q}$ for all $n\geq 0$.
\end{definition}

There is a good model structure, the Reedy structure, on cosimplicial spaces. Let $X$ be a cosimplicial space
so that $X^n$ is a simplicial set for $n\geq 0$. A morphism $f:X\rightarrow Y$ is a weak equivalence if each $f^n:X^n\rightarrow Y^n$ is a weak equivalence; that is, if $f$ is a level weak equivalence.
The maximal augmentation of a cosimplicial space is the simplicial set that equalizes $d^0,d^1:X^0\rightarrow X^1$. A map $f$ of cosimplicial spaces is called a cofibration if it is a level
cofibration (level monomorphism) and if it induces an isomorphism on the maximal augmentations.
The fibrations are all those morphisms with the right lifting property with respect to acyclic cofibrations. A proof that this
is a model category may be found in \cite[section X.5]{bousfield_homotopy_1972}.

\begin{example}\label{ex:cof}
    As examples of cofibrant objects, consider $\Delta$ and $\Delta[s]$. Indeed, $\Delta^0$ is a single point, and $\Delta^1$ is the $1$-simplex. The coface maps $d^0$ and $d^1$
    send the unique point of $\Delta^0$ to the vertices $1$ and $0$ respectively of $\Delta^1$. Therefore, the maximal augmentation is the empty simplicial complex.
    This also shows that $\Delta[s]\rightarrow\Delta$ is a cofibration.
\end{example}

Let $X$ be a cosimplicial space. Then, define the $n$th matching object of $X$ to be
\begin{equation*}
        M^n X=\underleftarrow{\lim}_{\phi:\mathbf{n}\rightarrow\mathbf{k}} X^k,
\end{equation*}
where $\phi$ runs over all surjections $\mathbf{n}\rightarrow\mathbf{k}$ in $\Delta$.
There is a natural map $X^{n+1}\rightarrow M^n X$.

\begin{proposition}[{\cite[Section~X.4.5]{bousfield_homotopy_1972}}]\label{prop:fibrants}
    A morphism $f:X\rightarrow Y$ is a fibration if and only if the induced map
    \begin{equation*}
            X^{n+1}\rightarrow Y^{n+1}\times_{M^n Y} M^n X
    \end{equation*}
    is a fibration of simplicial sets for all $n\geq -1$.
\end{proposition}

The closed model structure on cosimplicial spaces is simplicial. That is, there is a functor
\begin{equation*}
    \Map:\csSets^{\op}\times\csSets\rightarrow\sSets
\end{equation*}
defined by
\begin{equation*}
        \Map(X,Y):n\mapsto\Hom(X\times\Delta^n,Y).
\end{equation*}
The space $\Map(X,Y)$ is called the function complex from $X$ to $Y$.
Similarly, if $X$ and $Y$ are cosimplicial pointed spaces, then there is a pointed function complex
\begin{equation*}
    \Map_{\ast}(X,Y)\in\sSets_{\ast}
\end{equation*}
defined by
\begin{equation*}
    \Map_{\ast}(X,Y)_n=\Hom(X\wedge \Delta^{n}_+,Y).
\end{equation*}

\begin{proposition}[{\cite[Section~X.5]{bousfield_homotopy_1972}}]\label{prop:sm7}
The simplicial model category axiom SM7 is satisfied: if $A\rightarrow B$ is a cofibration of cosimplicial spaces and if $X\rightarrow Y$ is a fibration of cosimplicial spaces, then
\begin{equation*} 
        \Map(B,X)\rightarrow \Map(A,X)\times_{\Map(A,Y)}\Map(B,Y)
\end{equation*}
is a fibration.
\end{proposition}

There is a functor $\csSets\rightarrow\csSets_\ast$ defined by
\begin{equation*}
    X=(n\mapsto X^n)\mapsto X_+=(n\mapsto X^{n}_+),
\end{equation*}
where $X^{n}_+$ is the space $X^n$ with a disjoint basepoint attached.

\begin{definition}
  For each integer $n\geq 2$, there is a functor
  \begin{equation*}
    \pi_n:\csSets_\ast\rightarrow\cAb,
  \end{equation*}
  where $\cAb$ is the category of cosimplicial abelian groups,
  defined by
  \begin{equation*}
    \pi_n(X)^m=\pi_n(X^m).
  \end{equation*}
  There are also functors, defined by the same equation,
  \begin{align*}
    \pi_1:\csSets_\ast  &\rightarrow    \cGroups\\
    \pi_0:\csSets_\ast  &\rightarrow    \cSets_\ast,
  \end{align*}
  where $\cSets_\ast$ is the category of cosimplicial pointed sets and $\cGroups$ is the category of cosimplicial groups.
\end{definition}

\begin{definition}
  Let $A$ be a cosimplicial abelian group, cosimplicial group, or cosimplicial pointed set.
  A pointed cosimplicial space $X$ is called a $K(A,n)$-cosimplicial space if $\pi_n X\iso A$, while $\pi_m X\iso\ast$ for $m\neq n$. 
\end{definition}

\section{Spectral Sequences for Cosimplicial Spaces}\label{sec:ss}
Let $X$ be a pointed cosimplicial space. Define pointed simplicial complexes
\begin{align*} 
        \Tot_{\infty} X &= \Map_\ast(\Delta_+,X),\\
        \Tot_{s} X      &= \Map_\ast(\Delta[s]_+,X).
\end{align*}
By axiom SM7 and Example~\ref{ex:cof}, if $X$ is fibrant, then $\Tot_{s} X\rightarrow \Tot_{s-1} X$ gives a tower of pointed fibrations. The inverse limit of this tower is $\Tot_{\infty} X$ when $X$ is fibrant.

\begin{definition}\label{ss:total}
    For an arbitrary cosimplicial pointed space $X$ let $X\rightarrow\HH_c X$ be a pointed fibrant resolution. Then, the \df{total space spectral sequence} of $X$, $^\T\Eoh_1 X$, is defined
    to be the spectral sequence of the tower $\Tot_{*}\HH_c X$: 
    \begin{equation*}
        ^\T\Eoh_1^{s,t}X=\Eoh_{1}^{s,t}\{\Tot_{\ast}\HH_c X\}\rightharpoonup \pi_{t-s}\Tot_{\infty}\HH_c X.
    \end{equation*}
\end{definition}
    
The fiber $F_s$ of $\Tot_s \HH_c X\rightarrow\Tot_{s-1}\HH_c X$ and the homotopy groups of the fiber $F_s$ are identified in \cite[Proposition~X.6.3]{bousfield_homotopy_1972}
(see also Section~\ref{sec:differentials}):
\begin{equation*}
        F_s\simeq\Map_{*}(S^s,NX^s),
\end{equation*}
where $NX^s$ is the fiber of the fibration $\HH_c X^s\rightarrow M^{s-1}\HH_c X$; see Proposition \ref{prop:fibrants}. Moreover, 
\begin{equation*}
        \pi_i NX^s\iso \pi_i \HH_c X^s\cap\ker s^0\cap\cdots\cap\ker s^{s-1},
\end{equation*}
where the maps $s^i$ are the cosimplicial degeneracies:
\begin{equation*}
    s^i:\HH_c X^s\rightarrow\HH_c X^{s-1}.
\end{equation*}
Therefore, there is a natural identification
\begin{align*}
        ^\T\Eoh_1^{s,t}X\iso\pi_t \HH_c X^s\cap\ker s^0\cap\cdots\cap\ker s^{s-1} &\qquad\qquad t\geq s\geq 0.
\end{align*}
Note that for $t\geq 2$,
$^\T\Eoh_1^{s,t}X$ is the $s$-degree of the normalized cochain complex $N^*\pi_t\HH_c X$ associated to $\pi_t\HH_c X$.
It is tedious but not hard to check that under this identification, the differential $d_1$ of the spectral
sequence is chain homotopic to the differential of the normalized cochain complex. Therefore, there are natural identifications
\begin{align}\label{eq:totalsse2}
        ^\T\Eoh_2^{s,t}X\iso \Hoh^s(N^*\pi_t\HH_c)\iso \Hoh^s(C^* \pi_t \HH_c X)\iso \Hoh^s(C^* \pi_t X)  &\qquad    t\geq 2,\; t\geq s\geq 0,
\end{align}
where $C^*\pi_t X$ denotes the unnormalized cochain complex associated to the cosimplicial abelian group $\pi_t X$.

It is not hard to extend these identifications for
$t=0$ and $s=0$, and for $t=1$ and $s=0,1$. For a detailed discussion, see \cite[Section~X.7]{bousfield_homotopy_1972}.

Define, for a cosimplicial abelian group $A$, the $s$th cohomotopy group for $s\geq 0$ as
\begin{equation*}
        \pi^s A=\Hoh^s(C^* A).
\end{equation*}
If $G$ is a cosimplicial pointed set or cosimplicial group,
then define $\pi^0 G$ as the equalizer of $\partial^0,\partial^1:G^0\rightrightarrows G^1$. This is a group if $G$ is.

Similarly, define a pointed cohomotopy set $\pi^1 G$ for $G$ a cosimplicial group as follows.
Let
\begin{equation*}
        Z^1 G=\{g\in G^1: (\partial^0 g)(\partial^1 g)^{-1}(\partial^2 g)=1\}.
\end{equation*}
There is an action of $G^0\times Z^1 G\rightarrow Z^1 G$ given by
\begin{equation*}
        (g_0,g_1)\mapsto (\partial^1 g_0)g_1(\partial^0 g_0)^{-1}.
\end{equation*}
The set $Z^1 G$ is pointed by the element $1\in G^1$.
Let $\pi^1 G$ be the quotient set of $Z^1 G$ by this action, pointed by the orbit of $1\in Z^1 G$. Then, by \cite[Paragraph~X.7.2]{bousfield_homotopy_1972}, there are natural identifications
\begin{align}\label{eq:2tss}
    ^\T\Eoh_2^{s,t}X\iso
    \begin{cases}
        \pi^s\pi_t X    &   \text{if $t\geq s\geq 0$,}\\
        0               &   \text{otherwise.}
    \end{cases}
\end{align}

There is another spectral sequence for cosimplicial spaces, which is useful when $X$ is level Kan. This is the homotopy limit spectral sequence, which, in fact,
exists in much greater generality. See \cite[Chapter~XI]{bousfield_homotopy_1972}. The main tool is a functor from cosimplicial spaces to cosimplicial spaces called
$\Pi$. Let $N\Delta$ be the nerve of the category $\Delta$. Then an element of $N\Delta_n$ is a simplex $i_{\ast}:i_0\rightarrow\cdots\rightarrow i_n$, where
each $i_k$ is non-negative integer, and the
the arrows are order-preserving maps of the ordered sets associated to $i_k$: $\{0,\cdots,k\}$. For $X$ an arbitrary cosimplicial space,
let $\Pi X$ denote the space whose $n$th level is
\begin{equation*}
        \Pi^{n} X=\prod_{i_*\in N\Delta_n} X^{i_n}.
\end{equation*}
So, the $n$th level of $\Pi X$ is the product over all compositions $i_0\rightarrow i_n$ of $X^{i_n}$.
The face map $\partial^j$ for $j<n$ composed with projection onto $i_*$ is projection onto $\partial_j(i_*)$ followed by the identity. The face map $\partial^n$ composed with projection onto
$i_*$ is
projection onto $\partial_n(i_*)$ followed by $X(i_{n-1}\rightarrow i_n)$.
Similarly, the degeneracy $s^j$ followed by projection onto $i_*$ is projection onto $s_j(i_*)$ followed by the identity.

The important thing about the cosimplicial replacement functor $\Pi$ is that it takes level fibrations into cosimplicial fibrations and preserves weak equivalences.
See \cite[Proposition~X.5.3]{bousfield_homotopy_1972}.

\begin{definition}\label{ss:homotopy}
    Let $X$ be a cosimplicial pointed space. Let $\Ex^{\infty}X$ denote the cosimplicial pointed space obtained from $X$ by applying the $\Ex^{\infty}$-functor to each level.
    Then, $\Pi\Ex^{\infty}X$ is fibrant. Define the \df{homotopy limit spectral sequence} of $X$, $^{\HL}\Eoh_1 X$, to be
    \begin{equation*}
        ^{\HL}\Eoh_1^{s,t}X=\Eoh_1^{s,t}\{\Tot_{\ast}\Pi\Ex^{\infty}X\}\rightharpoonup\pi_{t-s}\Tot_{\infty}\Pi\Ex^{\infty} X.
    \end{equation*}
\end{definition}

The space $\Tot_{\infty}\Pi\Ex^{\infty} X$ is called the \df{homotopy limit} of $X$, and will be written as $\holim_{\Delta} X$.

\begin{lemma}
    If $X$ is a cosimplicial pointed space that is level Kan, then the natural morphism $X\rightarrow\Ex^{\infty} X$ induces an isomorphism of spectral sequences
    \begin{equation*}
        \Eoh_1^{s,t}\{\Tot_{\ast}\Pi X\}\riso\Eoh_1^{s,t}\{\Tot_{\ast}\Pi\Ex^{\infty} X\}.
    \end{equation*}
    This morphism is natural in morphisms of cosimplicial pointed level Kan spaces.
\end{lemma}

Let $X$ be an arbitrary pointed cosimplicial space.
The functor $\Pi$ can be defined on cosimplicial objects in any category with finite products. In particular, on pointed sets, groups, and abelian groups. There are natural isomorphisms,
of cosimplicial pointed sets for $n=0$, cosimplicial groups for $n=1$, and cosimplicial abelian groups for $n>1$,
\begin{equation*}
        \pi_n\Pi X\simeq \Pi \pi_n X,
\end{equation*}
where $\pi_n X$ is the cosimplicial object obtained by evaluating $\pi_n$ at each cosimplicial level.

For $X$ a cosimplicial object in a category with finite products, there is an natural morphism $X\rightarrow\Pi X$. The maps
\begin{equation*}
        X^n\rightarrow\prod_{i_*\in N\Delta_n} X^{i_n}
\end{equation*}
are described as follows.
The simplex $i_*$ determines a morphism $\left(0\rightarrow 1\rightarrow\cdots\rightarrow n\right)\rightarrow i_n$,
by taking the images of $0$ from each $i_i$, $0\leq i<n$. This induces the map $X^n$ to the product, and it extends to a cosimplicial map $X\rightarrow\Pi X$.

\begin{proposition}[{\cite[Paragraph~XI.7.3]{bousfield_homotopy_1972}}]
    The canonical map $C^\ast\pi_n X\rightarrow C^\ast\Pi\pi_n X$ is a quasi-isomorphism.
\end{proposition}

\begin{proposition}[{\cite[Paragraph~XI.7.5]{bousfield_homotopy_1972}}]
    If $X$ is a fibrant cosimplicial pointed space, then the natural morphism
    \begin{equation*}
        ^\T\Eoh_1 X\rightarrow ^{\HL}\Eoh_1 X
    \end{equation*}
    of spectral sequences is an isomorphism.
\end{proposition}

There is also the Postnikov tower of a cosimplicial space.

\begin{definition}\label{ss:postnikov}
    Let $X$ be a level-fibrant pointed cosimplicial space. Denote by $X(n)$ the level-wise application of the coskeleton functor. Then each $X\rightarrow X(n)$ and
    $X(m)\rightarrow X(n)$, $m\geq n$ is a level fibration. Therefore, $\Pi X(m)\rightarrow\Pi X(n)$ is a fibration for $m\geq n$. The spectral sequence of this tower is called 
    the \df{Postnikov spectral sequence} for $X$:
    \begin{equation*}
        ^{\Po}\Eoh_2^{s,t}X=\tilde{\Eoh}_2^{s,t}\{\Tot_{\infty}\Pi X(\ast)\}\iso\pi_{t-s}\Tot_{\infty}G(t)\rightharpoonup \pi_{t-s}\holim_{\Delta} X,
    \end{equation*}
    where $G(t)$ is the fiber of $\Pi X(t)\rightarrow\Pi X(t-1)$. This fiber is a fibrant resolution of a cosimplicial $K(\pi_t X,t)$-space.
\end{definition}
By \cite[Paragraphs~XI.7.2-3]{bousfield_homotopy_1972}, there are natural isomorphisms
\begin{equation*}
    \pi_{t-s} \Tot_{\infty}G(t)\simeq\pi^s\pi_t X
\end{equation*}
for $t\geq s\geq 0$.
Thus,
\begin{equation*}
    ^{\HL}\Eoh_2^{s,t}X\iso{^{\Po}\Eoh_2^{s,t}X}.
\end{equation*}
In fact, this isomorphism comes from an isomorphism of spectral sequences.

\begin{theorem}\label{thm:isomorphism}
    Let $X$ be level Kan. Then, there is a natural isomorphism $\phi$ of spectral sequences
    \begin{equation*}
        \phi:^{\HL}\Eoh_2 X\rightarrow ^{\Po}\Eoh_2 X
    \end{equation*}
    from the homotopy limit spectral sequence beginning with the $\Eoh_2$-page to the Postnikov tower spectral sequence.
    \begin{proof}
        Recall that to create an isomorphism of spectral sequences that come from exact couples
        \begin{center}
            \begin{tikzpicture}[description/.style={fill=white,inner sep=2pt}]
                    \matrix (m) [matrix of math nodes, row sep=3em,
                    column sep=2.5em, text height=1.5ex, text depth=0.25ex]
                    {     ^I\Doh_2 & & ^I\Doh_2\\
                          & ^I\Eoh_2, &\\};
                    \path[->,font=\scriptsize]
                    (m-1-1) edge node[auto] {$ i $} (m-1-3)
                    (m-1-3) edge node[below right] {$ j $} (m-2-2)
                    (m-2-2) edge node[below left] {$ k $} (m-1-1);
            \end{tikzpicture}
            \begin{tikzpicture}[description/.style={fill=white,inner sep=2pt}]
                    \matrix (m) [matrix of math nodes, row sep=3em,
                    column sep=2.5em, text height=1.5ex, text depth=0.25ex]
                    {     ^{II}\Doh_2 & & ^{II}\Doh_2\\
                    & ^{II}\Eoh_2 &\\};
                    \path[->,font=\scriptsize]
                    (m-1-1) edge node[auto] {$ i $} (m-1-3)
                    (m-1-3) edge node[below right] {$ j $} (m-2-2)
                    (m-2-2) edge node[below left] {$ k $} (m-1-1);
            \end{tikzpicture}
        \end{center}
        it suffices to create a morphism of exact couples that is an isomorphism just on the $\Eoh$-terms:
        \begin{equation*}
            \phi:^I\Eoh_2\riso{^{II}\Eoh_2}.
        \end{equation*}
        Indeed, this is enough to guarantee that the morphism induces an isomorphism on $\Hoh(\Eoh)$ and a morphism of the derived couples. So, it follows inductively that this is sufficient.

        Since $X$ is level Kan, there is a double tower of fibrations, of which a typical square is
        \begin{equation*}
                \begin{CD}
                        \Tot_{s+1}\Pi X(t-1)  @<<<  \Tot_{s+1}\Pi X(t)\\
                        @VVV                          @VVV\\
                        \Tot_{s}\Pi X(t-1)    @<<<  \Tot_{s}\Pi X(t).
                \end{CD}
        \end{equation*}
        These fit into a bigger diagram 
        \begin{equation}\label{dia:double}
                \begin{CD}
                        \Tot_{\infty}\Pi X(t-1)    @<<<  \Tot_{\infty}\Pi X(t)  @<<<  \Tot_{\infty}\Pi X\\
                        @VVV                              @VVV                          @VVV\\
                        \Tot_{s+1}\Pi X(t-1)      @<<<  \Tot_{s+1}\Pi X(t)    @<<<  \Tot_{s+1}\Pi X\\
                        @VVV                              @VVV                          @VVV\\
                        \Tot_s\Pi X(t-1)          @<<<  \Tot_s\Pi X(t)        @<<<  \Tot_s\Pi X.
                \end{CD}
        \end{equation}
        The horizontal inverse limits are $\Tot_{s}\Pi X$ and the vertical inverse limits are $\Tot_{\infty}\Pi X(t)$. Thus the homotopy limit spectral sequence comes the tower of fibrations
        at the horizontal limit, while the Postnikov spectral sequence comes from the tower of fibrations at the vertical limit. Define $F(s+1)$ to be the fiber of $\Tot_{s+1}\Pi X\rightarrow\Tot_s\Pi X$,
        and define $G(t)$ be the fiber of $\Tot_{\infty}\Pi X(t)\rightarrow\Tot_{\infty}\Pi X(t-1)$.

        First, construct a morphism
        \begin{equation*}
            ^{\HL}\Doh_2^{s,t}=\im(\pi_{t-s}\Tot_{s+1}\Pi X\rightarrow\pi_{t-s}\Tot_s\Pi X)\rightarrow{^\Po\Doh_2^{s,t}}=\pi_{t-s}\Tot_{\infty}\Pi X(t).
        \end{equation*}
        Let $[x]\in{^{\HL}\Doh_2^{s,t}}$ be represented by $x:S^{t-s}\rightarrow\Tot_{s+1}\Pi X\rightarrow\Tot_s\Pi X$. By adjunction, view this as
        \begin{equation*}
                x:\Delta[s]_+\wedge S^{t-s}\rightarrow\Delta[s+1]_+\wedge S^{t-s}\rightarrow\Pi X.
        \end{equation*}
        To create the morphism, one
        must ``lift'' this to a morphism $\phi(x):\Delta_+\wedge S^{t-s}\rightarrow\Pi X(t)$. The morphism $x$ consists of a compatible collection of morphisms
        \begin{equation*}
            x^{n,i_\ast}:\Delta^n[s]_+\wedge S^{t-s}\rightarrow\Delta^n[s+1]_+\wedge S^{t-s}\rightarrow X^{i_n},
        \end{equation*}
        one for each $i_\ast\in N\Delta_n$.
        The key point underlying the details below is that $X^{i_n}(t)$ is a Kan complex and also has trivial homotopy groups $\pi_k X^{i_n}(t)$ when $k> t$. At various points one needs to make choices
        to extend maps. These need to be compatible with the cosimplicial structure. At any given point, this will involve finitely many choices differing in simplicial degrees greater than $t$. Thus,
        it will always be possible to make the choices compatibly.

        Define $$\phi(x)^{0,i_\ast}:\Delta^0[r]_+\wedge S^{t-s}\rightarrow X^{i_0}(t)$$ for all $r$ and $i_\ast\in N\Delta_0$
        as the composition
        \begin{equation*}
            \Delta^0_+\wedge S^{t-s}=\Delta^0[s+1]_+\wedge S^{t-s}\xrightarrow{x^{0,i_\ast}} X^{i_0}\rightarrow X^{i_0}(t).
        \end{equation*}
        Since the coskeleton map is a functor, this definition is functorial.

        Now, suppose that
        \begin{equation*}
            \phi(x)^{h,i_*}:\Delta^{h}_+\wedge S^{t-s}\rightarrow X^{i_h}(t),
        \end{equation*}
        is defined for $0\leq h\leq k$ and all $i_*\in N\Delta_{h}$, compatible with all coface and codegeneracy maps in this range.
        In other words, suppose that we have defined compatible maps
        \begin{equation*}
            \phi(x)^{h,*}:\Delta^h_+\wedge S^{t-s}\rightarrow\Pi^h X(t)
        \end{equation*}
        for $0\leq h\leq k$.
        If $i_\ast\in N\Delta_{k+1}$ is degenerate, then define
        \begin{equation*}
            \phi(x)^{k+1,i_\ast}:\Delta^{k+1}_+\wedge S^{t-s}\rightarrow X^{i_{k+1}}(t)
        \end{equation*}
        by forcing the diagram
        \begin{equation*}
            \begin{CD}
                \Delta^{k+1}_+\wedge S^{t-s}    @>\phi(x)^{k+1,i_{\ast}} >>    X^{i_{k+1}}(t)\\
                @V s^j VV                                                   @|\\
                \Delta^k_+\wedge S^{t-s}        @>\phi(x)^{k,i'_\ast} >>    X^{i_{k+1}}(t)
            \end{CD}
        \end{equation*}
        to be commutative for every $j$ such that $s_j(i'_\ast)=i_\ast$ for some $i'_\ast\in N\Delta_k$. Making these simultaneously commutative
        for all choices of $j$ is possible because of the simplicial relations $s_i\circ s_j=s_{j+1}\circ s_i$ for $i\leq j$.

        If $i_\ast\in N\Delta_{k+1}$ is not degenerate, then the cosimplicial structure requires that the diagrams
        \begin{equation*}
            \begin{CD}
                \Delta^k_+\wedge S^{t-s}    @>\phi(x)^{k,\partial_j(i_\ast)} >>     X^{i_{k+1}}(t)\\
                @V\partial^j VV                                                     @|\\
                \Delta^{k+1}_+\wedge S^{t-s}@>\phi(x)^{k+1,i_\ast} >>               X^{i_{k+1}}(t)
            \end{CD}
        \end{equation*}
        for $0\leq j< k$ and
        \begin{equation*}
            \begin{CD}
                \Delta^k_+\wedge S^{t-s}    @>\phi(x)^{k,\partial_k(i_\ast)} >>      X^{i_k}(t)\\
                @V\partial^k VV                                                     @V X(i_{k+1}\rightarrow i_k)  VV\\
                \Delta^{k+1}_+\wedge S^{t-s}@>\phi(x)^{k+1,i_\ast} >>               X^{i_{k+1}}(t)
            \end{CD}
        \end{equation*}
        be commutative. In other words, the map $\phi(x)^{k+1,i_\ast}$ is already determined on $\partial\Delta^{k+1}_+\wedge S^{t-s}$.

        Thus, to make the inductive step, one must fill in the dashed line so that the diagram
        \begin{center}
            \begin{tikzpicture}[description/.style={fill=white,inner sep=2pt}]
                    \matrix (m) [matrix of math nodes, row sep=3em,
                    column sep=2.5em, text height=1.5ex, text depth=0.25ex]
                    {   \Delta^{k+1}[s+1]_+\wedge S^{t-s}       &   X^{i_{k+1}} & \\
                        \Delta^{k+1}_+\wedge S^{t-s}            &               &   X^{i_{k+1}}(t)    \\
                        \partial\Delta^{k+1}_+\wedge S^{t-s}    &               &\\};
                    \path[->,font=\scriptsize]
                    (m-1-1) edge node[auto] {$ x^{k+1,i_\ast} $} (m-1-2)
                    (m-1-2) edge node[auto] {$  $} (m-2-3)
                    (m-3-1) edge node[auto] {$  $} (m-2-3);
                    \path[right hook->,font=\scriptsize]
                    (m-1-1) edge node[auto] {$  $} (m-2-1)
                    (m-3-1) edge node[auto] {$  $} (m-2-1);
                    \path[dashed,->]
                    (m-2-1) edge node[auto] {$ \phi(x)^{k+1,i_\ast}  $} (m-2-3);
            \end{tikzpicture}
        \end{center}
        is commutative. If $s\geq k$, then $\Delta^{k+1}[s+1]_+=\Delta^{k+1}_+$, so there is nothing to do. If $s<k$, then $\Delta^{k+1}[s+1]_+\wedge S^{t-s}\subseteq\partial\Delta^{k+1}_+\wedge S^{t-s}$,
        and the outer square commutes by induction.
        Therefore, it suffices for the dashed arrow to commute in the bottom triangle. As $s<k$, the arrow $\partial\Delta^{k+1}_+\wedge S^{t-s}\rightarrow X^{i_{k+1}}$
        corresponds to a map $S^{k+t-s+1}\rightarrow X^{i_{k+1}}(t)$. But, $k+t-s+1> t$. Hence, $\pi_{k+t-s+1}X^{i_{k+1}}(t)=0$. As $X(t)$ is a Kan complex, such a fill exists
        (see \cite[page~35]{goerss_simplicial_1999}).

        Choose a fill for all choices of $i_\ast$. This completes the induction, giving
        \begin{equation*}
            \phi(x)^\ast:\Delta^{\ast}_+\wedge S^{t-s}\rightarrow\Pi^\ast X(t),
        \end{equation*}
        for $0\leq\ast\leq k+1$. The process outlined earlier in the proof gives a base case. So, induction provides the desired map
        \begin{equation*}
            \phi(x):\Delta_+\wedge S^{t-s}\rightarrow\Pi^\ast X(t).
        \end{equation*}

        If $y$ is another morphism $S^{t-s}\rightarrow\Tot_{s+1}\Pi X$ representing the class $[x]$,
        then it is straightforward, using a similar inductive argument, to lift a homotopy between $x$ and $y$ to a homotopy between $\phi(x)$ and $\phi(y)$
        so that the map is well-defined on the level of homotopy groups. Thus, the construction above determines a  well-defined the map $\phi:^{\HL}\Doh_2^{s,t}\rightarrow {^\Po\Doh_2^{s,t}}$. 

        It is easy to see that $\phi$ commutes with $i$. Indeed, on $^{\HL}\Doh_2$, $i$ is given by restriction from $\Delta[s]$ to $\Delta[s-1]$. On $^{\Po}\Doh_2$,
        $i$ is given by mapping $X(t)\rightarrow X(t-1)$. Now, composing $\phi(x)$ with $\Pi X(t)$ gives a map that solves the lifting problem to define $\phi(i(x))$:
        \begin{equation*}
            \begin{CD}
                \Delta[s-1]_+\wedge S^{t-s}     @>i(x) >>   \Pi X\\
                @VVV                                        @|\\
                \Delta[s]_+\wedge S^{t-s}       @>x >>      \Pi X\\
                @VVV                                        @|\\
                \Delta[s+1]_+\wedge S^{t-s}     @>>>        \Pi X\\
                @VVV                                        @VVV\\
                \Delta_+\wedge S^{t-s}          @>\phi(x) >>    \Pi X(t)\\
                @|                                         @VVV\\
                \Delta_+\wedge S^{t-s}          @>\phi(i(x)) >> \Pi X(t-1).
            \end{CD}
        \end{equation*}
        Therefore, $\phi$ commutes with $i$.

        Extend $\phi$ to $^{\HL}\Eoh_2$ by lifting
        \begin{equation*}
                x:S^{t-s}\rightarrow F(s)
        \end{equation*}
        to
        \begin{equation*}
                x:S^{t-s}\rightarrow\Tot_{s+1}\Pi X
        \end{equation*}
        and applying $\phi$ (recall the description in Section~\ref{sub:towers}  of $\Eoh_2$). In other words, apply $k$ and then $\phi$. This defines an element
        of $^{\Po}\Doh_2^{s,t}$, and, by construction, $\phi(x)$ actually factors through the fiber $G(t)$. Indeed, since
        \begin{equation*}
                x:\Delta[s+1]\wedge S^{t-s}\rightarrow\Pi X
        \end{equation*}
        comes from the fiber $F_s$, it restricts to the trivial map on
        \begin{equation*}
                \Delta[s-1]\wedge S^{t-s}.
        \end{equation*}
        Therefore, $x$ is trivial on the $t-1$-skeleton. Extending this to a map
        \begin{equation*}
            x:\Delta\wedge S^{t-s}\rightarrow\Pi X(t)
        \end{equation*}
        as above does not change this, so that when one composes with $\Pi X(t)\rightarrow\Pi X(t-1)$, the map is homotopic to the constant map on the basepoint.
        To check that this determines a well-defined map on $\Eoh_2$-terms, one must check that if
        \begin{equation*}
            x=y+j(z),
        \end{equation*}
        where $z\in\ker(\pi_{t-s+1}X_{s-1}\rightarrow\pi_{t-s+1}X_{s-2})$ and
        \begin{equation*}
            x,y\in\ker(\pi_{t-s}F(s)\rightarrow\pi_{t-s}X_s/i(\pi_{t-s}X_{s+1})),
        \end{equation*}
        then $\phi(x)=\phi(z)$. But, since
        \begin{equation*}
            \phi(j(z))=\phi(k(j(z)))=\phi(0),
        \end{equation*}
        it follows that this definition of $\phi$ on $^{\HL}\Eoh_2$ is indeed well-defined.

        It remains only to check that $\phi$ respects $j$. Again, let $[x]\in\pi_{t-s}\Tot_s\Pi X$ be represented by
        \begin{equation*}
                x:\Delta[s+1]_+\wedge S^{t-s}\rightarrow\Pi X.
        \end{equation*}
        Recall that $j([x])$ is obtained by
        the map
        \begin{equation*}
                \pi_{t-s}\Tot_{s+1}\Pi X\rightarrow\pi_{t-s-1} F(s+2),
        \end{equation*}
        given by lifting a horn
        \begin{equation*}
            \Delta[s+1]_+\wedge\Lambda_{i,+}^{s-t}\rightarrow\Delta[s+1]_+\wedge S^{s-t}\rightarrow \Pi X
        \end{equation*}
        to
        \begin{equation*}
                \Delta[s+2]_+\wedge\Delta^{s-t}_+\rightarrow\Pi X
        \end{equation*}
        and then restricting to the $i$th face to obtain
        \begin{equation*}
                \Delta[s+2]_+\wedge S^{t-s-1}\rightarrow\Pi X,
        \end{equation*}
        which by adjunction is in $\pi_{t-s-1} F(s+2)$.
        Lift
        \begin{equation*}
                \Delta[s+2]_+\wedge\Delta^{s-t}_+\rightarrow\Pi X
        \end{equation*}
        as above to a map
        \begin{equation*}
                \Delta_+\wedge\Delta^{s-t}_+\rightarrow\Pi X(t+1).
        \end{equation*}
        By construction, it maps down to
        \begin{equation*}
                \phi(x):\Delta_+\wedge S^{t-s}\rightarrow\Pi X(t),
        \end{equation*}
        while the restriction to $\Delta_+\wedge S^{t-s-1}$ is in the fiber $G(t+1)$. Therefore, $\phi$ commutes with $j$.

        To prove that $\phi$ is injective on $\Eoh_2$, let $x\in{^{\HL}\Eoh_2^{s,t}}$ be represented by
        \begin{equation*}
            x:\Delta[s+1]_+\wedge S^{t-s}\rightarrow\Pi X,
        \end{equation*}
        and let
        \begin{equation*}
            \phi(x):\Delta_+\wedge S^{t-s}\rightarrow\Pi X(t)
        \end{equation*}
        factor through $G(t)$. Suppose that $\phi(x)$ is homotopic to the constant map in $G(t)$, and let
        \begin{equation*}
            y:\Delta^1_+\wedge\Delta_+\wedge S^{t-s}\rightarrow G(t)\rightarrow\Pi X(t)
        \end{equation*}
        be such a homotopy. It is possible to find an extension
        \begin{center}
            \begin{tikzpicture}[description/.style={fill=white,inner sep=2pt}]
                    \matrix (m) [matrix of math nodes, row sep=3em,
                    column sep=2.5em, text height=1.5ex, text depth=0.25ex]
                    {   \Delta^1_+\wedge\Delta[s]_+\wedge S^{t-s}       &   \Pi X  \\
                        \Delta^1_+\wedge\Delta_+\wedge S^{t-s}          &   \Pi X(t)    \\};
                    \path[->,font=\scriptsize]
                    (m-1-1) edge node[auto] {$  $} (m-2-1)
                    (m-2-1) edge node[auto] {$ y $} (m-2-2)
                    (m-1-2) edge node[auto] {$  $} (m-2-2);
                    \path[dashed,->]
                    (m-1-1) edge node[auto] {$  $} (m-1-2);
            \end{tikzpicture}
        \end{center}
        giving a homotopy between $x$ and the constant map by definition of the coskeleton $\Pi X(t)$. Therefore, $\phi$ is injective on $\Eoh_2$.

        To show that $\phi$ is surjective on $\Eoh_2$, let
        \begin{equation*}
            y:S^{t-s}\rightarrow G(t)
        \end{equation*}
        be represented by
        \begin{equation*}
            \begin{CD}
                \Delta_+\wedge S^{t-s}      @>y >>      \Pi X(t)\\
                @VVV                                    @VVV\\
                \ast                        @>>>        \Pi X(t-1).
            \end{CD}
        \end{equation*}
        A lift $x$ for the diagram
        \begin{center}
            \begin{tikzpicture}[description/.style={fill=white,inner sep=2pt}]
                    \matrix (m) [matrix of math nodes, row sep=3em,
                    column sep=2.5em, text height=1.5ex, text depth=0.25ex]
                    {   \Delta^n[s-1]_+\wedge S^{t-s}       &   \ast  \\
                        \Delta^n[s]_+\wedge S^{t-s}          &        \\
                        \Delta^n[s+1]_+\wedge S^{t-s}         &   X^{i_n}\\
                        \Delta^n_+\wedge S^{t-s}            &     X^{i_n}(t)\\
                        \ast                                &   X^{i_n}(t-1)\\};
                    \path[->,font=\scriptsize]
                    (m-1-1) edge node[auto] {$  $} (m-1-2)
                    (m-1-1) edge node[auto] {$  $} (m-2-1)
                    (m-2-1) edge node[auto] {$  $} (m-3-1)
                    (m-3-1) edge node[auto] {$  $} (m-4-1)
                    (m-1-2) edge node[auto] {$  $} (m-3-2)
                    (m-3-2) edge node[auto] {$  $} (m-4-2)
                    (m-4-1) edge node[auto] {$ y $} (m-4-2)
                    (m-2-1) edge node[auto] {$  $} (m-3-2)
                    (m-4-1) edge node[auto] {$  $} (m-5-1)
                    (m-4-2) edge node[auto] {$  $} (m-5-2)
                    (m-5-1) edge node[auto] {$  $} (m-5-2);
                    \path[dashed,->,font=\scriptsize]
                    (m-3-1) edge node[auto] {$ x $} (m-3-2);
            \end{tikzpicture}
        \end{center}
        exists by definition for $n\leq s$. For $n>s$, the top quadrilateral means that the map
        \begin{equation*}
            \Delta^n[s]_+\wedge S^{t-s}\rightarrow X^{i_n}
        \end{equation*}
        is a diagram of $t-s$-spheres in $\Omega^{s}X^{i_n}$, and finding a lift $x$ is the same as saying that this class is $0$ in $\pi_{t-s}\Omega^s X^{i_n}$.
        But, the map $y$ in the diagram says that this class \emph{is} $0$ in $\pi_{t-s}\Omega^s X^{i_n}(t)$. Once again, by definition of the coskeleton
        functor, it follows that it is already $0$ in $\pi_{t-s}\Omega^s X^{i_n}$. Now, arguing as above, one may inductively choose extensions in the diagram to create an element $x$ of $^{\HL}\Eoh_2^{s,t}$
        such that $\phi(x)=y$. Therefore, $\phi$ is also surjective on $\Eoh_2$.

    \end{proof}
\end{theorem}

\begin{remark}
    This theorem does not appear to be new. Indeed, it appears that Thomason was aware of it: see \cite[page 542, paragraph 3]{thomason_algebraic_1985}.
    However, we know of no reference for a proof. It in fact holds in the greater generality of presheaves of simplicial sets on small categories. 
\end{remark}

\begin{corollary}
    Suppose that $X$ is a level Kan cosimplicial pointed space. If one of the spectral sequences
    \begin{align*}
        ^{\HL}\Eoh_1 X                       &\rightharpoonup    \pi_{t-s}\holim_{\Delta} X\\
        ^{\Po}\Eoh_2 X                       &\rightharpoonup    \pi_{t-s}\holim_{\Delta} X.
    \end{align*}
    converges, then both do, and the filtrations $^{\HL}Q_s$ and $^{\Po}Q_s$ coincide on $\pi_\ast \holim_{\Delta} X$.
\end{corollary}

\section{The \v{C}ech Approximation}\label{sec:cech}
Let $C$ be a Grothendieck site with terminal object $U$. Denote by $\Pre(C)$ the category of presheaves on $C$, and write $\sPre(C)$
for the category of simplicial presheaves.

We use the following closed model category structure on simplicial presheaves, called the local model category structure.
The cofibrations are the pointwise cofibrations. Thus, $X\rightarrow Y$ is a cofibration if and only if
$X(V)\rightarrow Y(V)$ is a monomorphism for every object
$V$ of $C$. For an object $V$ of $C$, there is a site with terminal object $C_{/V}$. Each presheaf on $C$ restricts
to a presheaf on $C_{/V}$. For a simplicial presheaf $X$, an object $V$ of $C$, and a basepoint $x\in X(V)_0$, there are presheaves of homotopy groups $\pi_k^p(X|V,x)$:
\begin{equation*}
  (f:W\rightarrow V)\mapsto \pi_k(|X(W)|,f^{*}(x)),
\end{equation*}
where $|X(W)|$ denotes the geometric realization of the simplicial set $X(W)$. Let $\pi_k(X|V,x)$ be the associated homotopy sheaf.
Call $w:X\rightarrow Y$ a weak equivalence if it induces an isomorphism of homotopy sheaves
\begin{equation*}
  \pi_k(X|V,x)\riso\pi_k(Y|V,w(x))
\end{equation*}
for all choices of $V$, all basepoints $x$ of $X(V)$, and all $k\geq 0$. 
The fibrations are all maps having the right lifting property with respect to all cofibrations that are simultaneously weak equivalences (the acyclic cofibrations). 
That this is a simplicial closed model category is due to Joyal; for a proof, see Jardine's paper \cite{jardine_simplicial_1987}.
Refer to these classes of morphisms more specifically as \df{global fibrations}, \df{global cofibrations}, and \df{local weak equivalences}.

\begin{theorem}[Dugger-Hollander-Isaksen {\cite{dugger_hypercovers_2004}}]\label{thm:dhi}
    If $X$ is a simplicial presheaf, and if $\mathcal{V}^\bullet\rightarrow V$ is a hypercover in $C$, then let $X_{\mathcal{V}^\bullet}$ denote the cosimplicial space associated to $\mathcal{V}^{\bullet}$.
    There is a canonical augmentation $X(V)\rightarrow X_{\mathcal{V}^{\bullet}}$.
    The simplicial presheaf $X$ is globally fibrant if and only if
    \begin{equation*}
        X(V)\rightarrow \holim_{\Delta}X_{\mathcal{V}^{\bullet}}
    \end{equation*}
    is a weak equivalence for all hypercovers $\mathcal{V}^{\bullet}\rightarrow V$ in $C$.
\end{theorem}

There are other types of morphisms, namely pointwise
weak equivalences and pointwise fibrations. A pointwise weak equivalence is a morphism $f:X\rightarrow Y$ such that $X(V)\rwe Y(V)$ is a weak equivalence of simplicial sets for all objects $V$ of $C$.
Two pointwise weak equivalent sheaves are local weak equivalent, and two local weak equivalent fibrant presheaves are pointwise weak equivalent.
A pointwise fibration is a morphism $f:X\rightarrow Y$ such
that every $f:X(V)\rightarrow Y(V)$ is a fibration of simplicial sets. Say that $X$ is pointwise Kan or pointwise fibrant if $X\rightarrow\ast$ is a pointwise fibration

Note that if a simplicial presheaf is pointed, then the homotopy presheaves and sheaves above may be defined globally.

Let $F$ be a functor from simplicial sets to simplicial sets such that $F(\emptyset)=\emptyset$, or from pointed simplicial sets to pointed simplicial sets
such that $F(\ast)=\ast$. If $X$ is a simplicial presheaf,
denote by $FX$ the pointwise application of $F$ to $X$, so that $(FX)(V)=F(X(V))$ for all $V$. For instance, below, $\cosk_n X$ will be the pointwise $n$-coskeleton of $X$. If
$F$ preserves weak equivalences of simplicial sets, then it preserves local weak equivalences of simplicial presheaves. This is the case, for instance, for the coskeleta functors and for the $\Ex$ functor.
In particular, one may always replace $X$ with the pointwise weakly equivalent presheaf $\Ex^{\infty} X$, which is pointwise Kan.

\begin{definition}\label{def:hypercohomology}
    If $X$ is a presheaf of pointed spaces, and if $X\rightarrow\HH X$ is a pointed fibrant resolution of $X$, then $\HH X$ is called the \df{hypercohomology space} of $X$. Hypercohomology sets, groups, and
    abelian groups are defined by
    \begin{equation*}
        \HH^s(U,X)=\pi_s\Gamma(U,\HH X).
    \end{equation*}
\end{definition}

\begin{definition}\label{ss:brown-gersten}
    Let $X$ be a presheaf of pointed simplicial sets on $C$, and let $X\rightarrow\HH X$ be a pointed fibrant resolution.
    Finally, let $X(n)$ be a pointed fibrant resolution of $\cosk_n \HH X$ so that $X(n)\rightarrow X(n-1)$ is a global fibration for all $n\geq 0$.
    Then, the \df{Brown-Gersten spectral sequence associated to $X$} is the $\tilde{\Eoh}_2$-spectral sequence (see Section~\ref{sub:towers}) associated to the tower of fibrations
    \begin{equation*}
        \Gamma(U,X(n))\rightarrow\Gamma(U,X(n-1)).
    \end{equation*}
    That is, define the Brown-Gersten spectral sequence as
    \begin{equation*}
        ^{\BG}\Eoh_2^{s,t}X=\tilde{\Eoh}_2^{s,t}\{\Gamma(U,X(\ast))\}\rightharpoonup\pi_{t-s}\lim_{\leftarrow}\Gamma(U,X(n)).
    \end{equation*}
    This is also called the \df{descent spectral sequence}.
\end{definition}

\begin{lemma}[{\cite{gillet_filtrationshigher_1999}}]
    Suppose that the site $C$ is locally of finite cohomological dimension.
    Then, the natural morphism
    \begin{equation*}
        \Gamma(U,\HH X)\rightarrow \lim_{\leftarrow}\Gamma(U,X(n))
    \end{equation*}
    is a weak equivalence for all objects $U$ whenever $X$ is locally weak equivalent to a product of a constant pointed space and a connected space.
\end{lemma}

In the situation of the lemma, write the Brown-Gersten spectral sequence as
\begin{equation*}
    ^{\BG}\Eoh_2 X\rightharpoonup\HH^{t-s}(U,X).
\end{equation*}
For instance, the lemma holds for presheaves of Quillen $K$-theory spaces on the \'{e}tale sites of schemes of finite \'{e}tale cohomological dimension.

The theorem of Dugger-Hollander-Isaksen, Theorem~\ref{thm:dhi}, says that the natural map
\begin{equation*}
    \Gamma(U,X(n))\rightarrow\Tot_{\infty}\Pi X(n)_{\mathcal{U}^{\bullet}}
\end{equation*}
is an isomorphism whenever $\mathcal{U}^{\bullet}\rightarrow U$ is a hypercover, and $X(n)_{\mathcal{U}^{\bullet}}$ denotes the cosimplicial space obtained by evaluating $X(n)$ at each level of
$\mathcal{U}^{\bullet}$.

\begin{definition}
    Let $\mathcal{U}^{\bullet}$ be a hypercover of an object $U$ of $C$. Let $X$ be a simplicial presheaf. Then, let $X_{\mathcal{U}^{\bullet}}$
    denote the cosimplicial space one gets by evaluating $X$ at $\mathcal{U}^{\bullet}$.
    The \df{\Cech hypercohomology} of $X$ on $\mathcal{U}^{\bullet}$ is defined as
    \begin{equation*}
        \check{\HH}^s(\mathcal{U}^{\bullet},X)=\pi_{s}\Tot_{\infty}\HH_{c} X_{\mathcal{U}^{\bullet}}.
    \end{equation*}
\end{definition}

\begin{definition}\label{ss:cech}
    The \df{\Cech hypercohomology spectral sequence} associated to $X$ and $\mathcal{U}^{\bullet}\rightarrow U$ is the total space spectral sequence associated to $\HH_c X_{\mathcal{U}^{\bullet}}$:
    \begin{equation*}
        ^{\mathcal{U}^{\bullet}}\Eoh_1 X=^\T\Eoh_1 X_{\mathcal{U}^{\bullet}}.
    \end{equation*}
\end{definition}
By Equation~\eqref{eq:totalsse2}, the $\Eoh_2$-terms for the corresponding spectral sequence are naturally identified with
\begin{equation*}
    \Eoh_2^{s,t}\simeq\check{\Hoh}^s(\mathcal{U}^{\bullet},\pi_t^p X),
\end{equation*}
as desired.

There is a natural morphism from \Cech hypercohomology to hypercohomology induced by a natural morphism
\begin{equation}\label{eq:hypermorphism}
    \Tot_{\infty}\HH_c X_{\mathcal{U}^{\bullet}}\rightarrow\Gamma(U,\HH X).
\end{equation}
This morphism is the composition of
\begin{equation*}
    \Tot_{\infty}\HH_c X_{\mathcal{U}^{\bullet}}\rightarrow\Tot_{\infty}\HH_c\Pi^\ast X_{\mathcal{U}^{\bullet}}=\Tot_{\infty}\Pi^\ast X_{\mathcal{U}^{\bullet}}\rightarrow\Tot_\infty \Pi^\ast (\HH X)_{\mathcal{U}^{\bullet}}
\end{equation*}
with the inverse of the local weak equivalence (Theorem~\ref{thm:dhi})
\begin{equation*}
    \Gamma(U,\HH X)\rwe \Tot_{\infty}\Pi^\ast (\HH X)_{\mathcal{U}^{\bullet}}.
\end{equation*}

\begin{theorem}\label{thm:morphism}
    Let $X$ be a pointed simplicial presheaf, and let $\mathcal{U}^{\bullet}\rightarrow U$ be a hypercover.
        There is a natural morphism
        \begin{equation*}
            ^{\mathcal{U}^{\bullet}}\Eoh_2 X\rightarrow {^{\BG}\Eoh_2 X}
        \end{equation*}
        from the \Cech hypercohomology spectral sequence to the Brown-Gersten spectral sequence, which on $\Eoh_2$-terms is the natural map
        \begin{equation*}
            \check{\Hoh}^s(\mathcal{U}^{\bullet},\pi_{t}^p X)\rightarrow\Hoh^s(U,\pi_{t} X),
        \end{equation*}
        and which respects the morphism on abutments of Equation~\eqref{eq:hypermorphism}.
        \begin{proof}
            One may assume, possibly by applying the $\Ex^\infty$ that $X$ is pointwise Kan. Therefore, $X_{\mathcal{U}^{\bullet}}$ is a level Kan cosimplicial pointed space.
            Then, the space $\Pi X_{\mathcal{U}^{\bullet}}$
            is fibrant, so that there is a natural morphism
            \begin{equation*}
                \HH_c X_{\mathcal{U}^{\bullet}}\rightarrow\Pi X_{\mathcal{U}^{\bullet}}
            \end{equation*}
            making the diagram
            \begin{center}
                \begin{tikzpicture}[description/.style={fill=white,inner sep=2pt}]
                        \matrix (m) [matrix of math nodes, row sep=3em,
                        column sep=2.5em, text height=1.5ex, text depth=0.25ex]
                        {     X_{\mathcal{U}^{\bullet}} & & \Pi X_{\mathcal{U}^{\bullet}}\\
                        & \HH_c X_{\mathcal{U}^{\bullet}} &\\};
                        \path[->,font=\scriptsize]
                        (m-1-1) edge node[auto] {$  $} (m-1-3)
                        (m-1-1) edge node[below right] {$  $} (m-2-2)
                        (m-2-2) edge node[below left] {$  $} (m-1-3);
                \end{tikzpicture}
            \end{center}
            commutative.
            This induces a morphism of total space spectral sequences
            \begin{equation}\label{eq:morph1}
                ^{\mathcal{U}^{\bullet}}\Eoh_2 X\rightarrow{^\T\Eoh_2\HH_c X_{\mathcal{U}^{\bullet}}}\rightarrow {^{\HL}\Eoh_2 X_{\mathcal{U}^{\bullet}}}
            \end{equation}
            There is a morphism
            \begin{equation}\label{eq:morph2}
                ^{\HL}\Eoh_2 X_{\mathcal{U}^{\bullet}}\rightarrow {^{\HL}\Eoh_2 (\HH X)_{\mathcal{U}^{\bullet}}}
            \end{equation}
            induced by $X\rightarrow\HH X$. 
            Theorem \ref{thm:isomorphism} furnishes an isomorphism
            \begin{equation}\label{eq:morph3}
                ^{\HL}\Eoh_2(\HH X)_{\mathcal{U}^{\bullet}}\riso {^{\Po}\Eoh_2(\HH X)_{\mathcal{U}^{\bullet}}}.
            \end{equation}
            Again, the fibrant replacement $\cosk_n \HH X\rightarrow X(n)$ induces a morphism
            \begin{equation}\label{eq:morph4}
                ^{\Po}\Eoh_2(\HH X)_{\mathcal{U}^{\bullet}}\rightarrow\tilde{\Eoh}_2^{s,t}\{\Tot_{\infty}\Pi X(n)_{\mathcal{U}^{\bullet}}\}.
            \end{equation}
            Finally, the result of Dugger-Hollander-Isaksen \cite{dugger_hypercovers_2004} says that the natural morphism
            \begin{equation}\label{eq:morph5}
                ^{\BG}\Eoh_2 X\rightarrow\tilde{\Eoh}_2^{s,t}\{\Tot_{\infty}\Pi X(n)_{\mathcal{U}^{\bullet}}\}
            \end{equation}
            is in fact an isomorphism.

            The theorem follows by taking the composition of the morphisms of Equations \eqref{eq:morph1}, \eqref{eq:morph2}, \eqref{eq:morph3}, and \eqref{eq:morph4},
            and the inverse of the morphism of Equation \eqref{eq:morph5}.
        \end{proof}
\end{theorem}

\section{Differentials in the Holim Spectral Sequence}\label{sec:differentials}

\begin{lemma}[{\cite[Proposition X.6.3.i]{bousfield_homotopy_1972}}]\label{lem:fiber}
    Take $X$ to be a fibrant pointed cosimplicial space.
    There are equivalences $F(s)\simeq \Map_*(S^s,NX^s)$, where $F(s)$ is the fiber of $\Tot_s X\rightarrow \Tot_{s-1} X$, and $NX^s$ is the fiber of the
    fibration $X^s\rightarrow M^{s-1} X$.
    \begin{proof}[Sketch of proof]
        Let $\beta:\Delta^n\rightarrow F(s)$.
        By adjunction, this is a morphism $\Delta[s]_+\wedge \Delta^n_+\xrightarrow{\beta} X$ such that the restriction to $\Delta[s-1]_+\wedge\Delta^n_+$
        factors through the base-point. In particular,
        the level $s$ picture is a map of pointed spaces $\Delta^s_+\wedge\Delta^n_+\rightarrow X^s$ such that the restriction to
        the $s-1$-skeleton on the left $\Delta^s[s-1]_+\wedge\Delta^n_+$ factors through the base-point of $X^s$. Therefore, this is a morphism $S^s\wedge\Delta^n_+\rightarrow X$.
        That is, $\beta$ determines an $n$-simplex of $\Map_*(S^s,X^s)$.
        However, the degeneracy $s^i\beta:\Delta^{s-1}[s]_+\wedge\Delta^n_+\rightarrow X^{s-1}$ factors through $\Delta^{s-1}[s-1]_+\wedge\Delta^n_+$ and so is trivial.
        Therefore the $n$-simplex of $\Map_*(S^s,X^s)$ actually lives in $\Map_*(S^s,NX^s)$.
        Conversely, suppose that $\gamma$ is an $n$-simplex of $\Map_*(S^s,NX^s)$. Again, by adjunction, view this
        as a map $S^s\wedge\Delta^n_+\rightarrow NX^s\rightarrow X^s$. Extend this as follows into a map $\Delta[s]_+\wedge\Delta^n_+\rightarrow X$. The lift of $\gamma$
        to $\Delta^s_+\wedge\Delta^n_+$
        is the $s$-level of this morphism. Since $\lambda$ factors
        through $*$ for any degeneracy on $\Delta^s$, let $\Delta^k[s]_+\wedge\Delta^n_+\rightarrow *\rightarrow X^k$ for $k<s$. If $k>s$, use the following diagram:
        \begin{equation}\label{eq:liftingdiagram}
                \begin{CD}
                        \Delta^s[s]_+\wedge\Delta^n_+    @>\beta >>                      X^s\\
                        @V\partial VV                                                @V\partial VV\\
                        \Delta^k[s]_+\wedge\Delta^n_+    @>\partial\circ\beta\circ s >>  X^k\\
                        @V s VV                                                      @V s VV\\
                        \Delta^s[s]_+\wedge\Delta^n_+    @>\beta >>                      X^s,
                \end{CD}
        \end{equation}
        where the face map $\partial$ (resp. the degeneracy map $s$) is some composition of face (resp. degeneracy) maps such that the vertical compositions are the identity.
        This defines the extension inductively on the faces of $\Delta^k[s]$.
        But, one need only define it up to faces, since $k>s$. This gives us an $n$-simplex of $\Tot_s X$. Clearly, these constructions are mutually inverse.
    \end{proof}
\end{lemma}

It is not hard, using this identification, to show that the homotopy groups of $F(s)$ are the groups of the normalized cochain complexes associated to $X$.
That is
\begin{equation*}
        \pi_{t-s} F(s)\simeq \pi_{t-s}\Map_*(S^s,NX^s)\simeq\pi_t X^s\cap \ker s^0\cap\cdots\cap\ker s^{s-1},
\end{equation*}
where $\ker s^i$ is $\ker(s^i:\pi_t X^s\rightarrow \pi_t X^{s-1})$.
See \cite[proposition X.6.3]{bousfield_homotopy_1972}.

The differentials in the homotopy limit spectral sequence are those of the spectral sequence associated to a tower of fibrations. In particular, they arise from trying to lift simplices to
higher and higher levels of the tower. The recipe is as follows: let
\begin{equation*}
    \cdots\rightarrow X_{s+1}\rightarrow X_s\rightarrow X_{s-1}\rightarrow\cdots
\end{equation*}
be a tower of fibrations, and let $\delta:S^t\rightarrow X_{s-1}$ be a homotopy class. Writing $S^t$ as the quotient of $\Delta^t$ by its boundary, we get a map $\Delta^t\rightarrow X_{s-1}$
that is trivial on the boundary $\partial\Delta^t$. We can obviously then map any horn $\Lambda^{t}_i$ to $X^s$ by sending it to the base-point. Then, there is a lifting problem
\begin{equation*}
    \begin{CD}
        \Lambda^t_i     @>>>    X_s\\
        @VVV                    @VVV\\
        \Delta^t        @>>>    X_{s-1}.
    \end{CD}
\end{equation*}
Since the map on the right is a fibration, we do get a lift to a map $\Delta^t\rightarrow X_s$. However, it may no longer map the boundary of $\Delta^t$ to the base-point of $X_s$.
Instead, the $i$th face determines a map $\Delta^{t-1}\rightarrow X_s$ that \emph{is} trivial on the boundary. The corresponding class in $\pi_{n-1}F(s)$ is the obstruction to lifting the homotopy
class of $\delta$ to a homotopy class on $X_s$.

Below, the tower of fibrations is the total space tower associated to a fibrant cosimplicial space (see Section~\ref{sec:cosimplicial}).

\begin{remark}\label{remark:holimd1}
    We briefly indicate how to show that $d_1$ is homotopic to the cochain differential on the normalized cochain complex for $\pi_t X$. Let $\beta$ this time represent a class in
    $\pi_{t-s}\Map_*(S^s,NX^s)$. This determines a corresponding map $S^s\wedge S^{t-s}\rightarrow X^s$. Extending this, one gets a map $\Delta[s]_+\wedge S^{t-s}\rightarrow X$. Lift this,
    using some horn $\Lambda^{t-s}_i$,
    to a map $\Delta[s+1]_+\wedge\Delta^{t-s}_+\rightarrow X$. On $\Delta[s]_+\wedge\Delta^{t-s}_+$ this agrees with $\beta$. Restricting to
    the $i$th face of $\Delta^{t-s}$, one gets a map $\Delta[s+1]_+\wedge S^{t-s-1}\rightarrow X$. Then, look at just $\Delta^{s+1}_+\wedge S^{t-s-1}\rightarrow X^{s+1}$, and check that one lands in the fiber.
    The reason that this agrees with the alternating sum $d=\Sigma (-1)^i d^i$ is that the simplices $\Delta^{s+1}_+\wedge\Delta^{t-s}_+\rightarrow X^{s+1}$ give a homotopy between $d(\beta)$ and $d_1(\beta)$.
\end{remark}

\begin{construction}\label{const:holimd2}
    Let $\delta:S^t\rightarrow X^s$ represent a class $[\delta]$ of $\pi^s \pi_t X$, where $t-s\geq 0$,
    and suppose that $\delta$ factors as $\delta:S^{t-s}\rightarrow\Map_*(S^s,NX^s)$. That is, suppose that $\delta$ represents
    $[\delta]$ in the normalized chain complex. Then, as in the proof of Lemma~\ref{lem:fiber}, extend $\delta$ to a map
    \begin{equation*}
            \delta':\Delta[s]_+\wedge\Delta^{t-s}_+\rightarrow X.
    \end{equation*}
    Now, the restriction of this map to $\Delta[s]_+\wedge\partial\Delta^{t-s}_+$ factors through the base-point.
    Choosing a horn $\Lambda^{t-s}_i\subset\partial\Delta^{t-s}$ and the map $\Delta[s+1]_+\wedge\left(\Lambda^{t-s}_i\right)\rightarrow *\rightarrow X$, one gets,
    using the fact that $\Tot_{s+1}X\rightarrow\Tot_s X$ is a fibration, a lift to a map
    \begin{equation*}
            \gamma:\Delta[s+1]_+\wedge\Delta^{t-s}_+\rightarrow X
    \end{equation*}
    such that the following diagram commutes
    \begin{equation*}
            \begin{CD}
                    \Delta[s+1]_+\wedge\Delta^{t-s}_+  @>\gamma>>      X\\
                    @AAA                                          @|\\
                    \Delta[s]_+\wedge\Delta^{t-s}_+    @>\delta^{'}>>  X.
            \end{CD}
    \end{equation*}
    The fact that $\delta$ represents a cohomotopy class implies, by Remark~\ref{remark:holimd1}, that the restriction of $\gamma$ to
    \begin{equation*}
            \Delta^{s+1}_+\wedge\partial_i\Delta^{t-s}_+\rightarrow X^{s+1}
    \end{equation*}
    is contractible in $\Map_*(S^{s+1},NX^{s+1})$. Therefore, one can replace $\gamma$ by a homotopic map $\gamma'$ such that the restriction to
    $\gamma':\Delta[s+1]_+\wedge\partial\Delta^{t-s}_+\rightarrow X$ factors through the base-point.
    Thus, one can repeat the process, using a trivial map $\Delta[s+2]_+\wedge\left(\Lambda^{t-s}_j\right)_+\rightarrow *\rightarrow X$, to get a lift $\epsilon:\Delta[s+2]_+\wedge\Delta^{t-s}_+\rightarrow X$,
    using some horn $\Lambda^{t-s}_j$.
    The differential $d_2([\delta])$ is the class of $\epsilon:\Delta^{s+2}_+\wedge\partial_j\Delta^{t-s}_+\rightarrow X^{s+2}$ in $\pi^{s+2}\pi_{t-s-1} X$.
\end{construction}

Note the following, which will be an aid to making the extensions described above. We claim that to extend $\Delta[s]_+\wedge\Delta^{t-s}_+\rightarrow X$ to $\Delta[s+1]_+\wedge\Delta^{t-s}_+$ it is sufficient
to describe the extension of $\Delta^{s+1}[s]_+\wedge\Delta^{t-s}_+\rightarrow X^{s+1}$ to $\Delta^{s+1}_+\wedge\Delta^{t-s}_+$. Indeed, make the same argument as in the
proof of Lemma~\ref{lem:fiber}, especially the argument using Diagram~\eqref{eq:liftingdiagram}.


\section{Description of $d_2$ for the \Cech Spectral Sequence}\label{sec:cechd2}

The descriptions of the differentials above translate into the following theorem in the setting of the \Cech approximation.

\begin{theorem}\label{thm:cechd2}
    Let $X$ be a presheaf of pointed simplicial sets on the site $C$, and let $\mathcal{U}^{\bullet}:\mathcal{V}_A\rightarrow\mathcal{U}_I\rightarrow U$ be a $1$-hypercover in $C$.
    Thus, for $\alpha\in A$, $V_{ij}^{\alpha}\rightarrow U_{ij}$ is a covering morphism, where $U_{ij}=U_i\times_U U_j$. Then, for $t>0$, the differentials
    \begin{equation*}
        d_2:\check{\Hoh}^0(\mathcal{U}^{\bullet},\pi_t X)\rightarrow\check{\Hoh}^2(\mathcal{U}^{\bullet},\pi_{t+1} X)
    \end{equation*}
    can be described as follows. An element $[\delta]$ of $\check{\Hoh}^0(\mathcal{U}^{\bullet},\pi_t X)$ is represented by
    a map
    \begin{equation*}
        \delta:S^t\rightarrow X_{\mathcal{U}^{\bullet}}^0=\prod_{i\in I} X(U_i)
    \end{equation*}
    such that $\delta_i$ is homotopic to $\delta_j$ on $V_{ij}^{\alpha}$ for every $i,j,\alpha$. Pick a specific based homotopy
    \begin{equation*}
        y_{ij}^{\alpha}:\delta_i\rightarrow\delta_j
    \end{equation*}
    on $V_{ij}^{\alpha}$. This data determines a map $\partial\Delta^2_+\wedge\Delta^t_+\rightarrow X_{\mathcal{U}^{\bullet}}^2$
    such that on each face of $\Delta^2$, the component in $X(V_{ijk}^{\alpha\beta\gamma})$ is one of the homotopies $y_{ij}^{\alpha}$, $y_{jk}^{\beta}$, or $y_{ki}^{\gamma}$. Then, let $\Lambda^{t}_n\subseteq\Delta^t$ be a horn,
    and let
    \begin{equation*}
        \Delta^2_+\wedge\Delta^t_+\rightarrow X_{\mathcal{U}^{\bullet}}^2
    \end{equation*}
    be a fill of the horn. Then, the restriction to
    \begin{equation*}
        \Delta^2_+\wedge\partial^n\Delta^t_+
    \end{equation*}
    is in the class of $d_2([\delta])$.
    \begin{proof}
        The proof follows immediately from Construction~\ref{const:holimd2} and Definition~\ref{ss:cech}.

    \end{proof}
\end{theorem}


\section{Divisibility Theorem}\label{sec:division}


In this section, the \Cech spectral sequence is used to study the differentials $d_2^{\alpha}$ in the Brown-Gersten spectral sequence of twisted algebraic $K$-theory.

To a class $\alpha\in\Hoh^2(U_{\et},\Gm)$, one associates a stack $\StProjA$ of $\alpha$-twisted coherent sheaves locally free and of finite rank.
This is a stack of exact categories. The pointwise $K$-theory of this stack will be written $\K^{\alpha}$, where
\begin{equation*}
    \K^{\alpha}(V)=\K^Q(\StProjA_V)=BQ(\StProjA_V),
\end{equation*}
for $V\rightarrow U$, and where $BQ(\StProjA_V)$ is the classifying space of Quillen's category $Q(\StProjA_V)$ \cite{quillen_higher_1973}.
The homotopy presheaves are
\begin{equation*}
    \K^{\alpha}_k:V\mapsto\pi_{k+1}\K^{\alpha}(V),
\end{equation*}
and the associated sheaves are $\mathcal{K}^{\alpha}_k$. The differentials of the Brown-Gersten spectral sequence associated to $\K^{\alpha}$
are written $d_r^{\alpha}$. The idea of an $\alpha$-twisted
sheaf was apparently created by Giraud \cite{giraud_cohomologie_1971}, and was revived in the thesis of \caldararu \cite{caldararu_derived_2000}.
For a modern geometric approach, see Lieblich~\cite{lieblich_twisted_2008}. In the context of $K$-theory, see~\cite{antieau_index_2009}.

Let $\K^{\alpha,\et}=\HH\K^{\alpha}$ be a fibrant replacement for $\K^{\et}$ in the local model category
of presheaves of spaces on $U_{\et}$. When $U$ is of finite cohomological dimension, the Brown-Gersten spectral sequence of $\K^{\alpha}$,
\begin{equation*}
    \Eoh_2^{s,t}=\Hoh^s(U_{\et},\pi_t\K^{\alpha})=\Hoh^s(U_{\et},\mathcal{K}_{t-1}^{\alpha})\Rightarrow\pi_{t-s}\K^{\alpha,\et}(U),
\end{equation*}
converges strongly to the homotopy of $\K^{\alpha,\et}$.

\begin{definition}
    In \cite{antieau_index_2009}, natural isomorphisms $\mathcal{K}^{\alpha}_k\riso\mathcal{K}_k$ are given. It follows that, for $U$ geometrically connected,
    $\Hoh^0(U_{\et},\mathcal{K}^{\alpha}_0)\riso\ZZ$.
    The \'etale index of $\alpha$, denoted by $eti(\alpha)$, is defined to be the smallest integer $n\in\Hoh^0(U_{\et},\mathcal{K}^{\alpha}_0)$ such that $d_k^{\alpha}(n)=0$ for all $k\geq 2$.
    In other words, $eti(\alpha)$ is the generator of
    \begin{equation*}
        \Eoh_{\infty}^{0,1}\subseteq\Hoh^0(U_{\et},\mathcal{K}^{\alpha}_0)=\Hoh^0(U_{\et},\pi_1\K^{\alpha}).
    \end{equation*}
    Equivalently, $eti(\alpha)$ is the positive generator of the image of the rank map $\K^{\alpha,\et}_0(U)\rightarrow\ZZ$.
    The period of $\alpha$, denoted by $per(\alpha)$, is the order of $\alpha$ in $\Hoh^2(U_{\et},\Gm)$. 
\end{definition}

\begin{remark}
    In general, if $X$ is a presheaf of spaces, then the differentials leaving $\Hoh^0(U,\pi_0 X)$ in the Brown-Gersten spectral sequence for $X$ are all zero.
    This is why we use $BQ(\StProjA)$ instead of $\Omega BQ(\StProjA)$ as our model for $K$-theory. Of course, one can stabilize the argument and work with presheaves of spectra.
    But, that argument requires another level of complexity which we wished to avoid.
\end{remark}

A key ingredient of our main theorem is the following lemma, which allows us to identify the critical lift in the description of $d_2$ in the \Cech spectral sequence for twisted $K$-theory.
First, some notation.

Let $E$ be an exact category, let $M_{i}$, $i=0,1,2$, be objects of $E$, and let $\theta_{ij}:M_i\rightarrow M_j$ be isomorphisms for $\theta_{01}$, $\theta_{12}$, and $\theta_{20}$.
Then, the $M_i$ determine loops in $BQE$, that is, elements of $\pi_1 BQE$. Recall here that the base-point of $BQE$ is the zero object of $E$.
The isomorphisms $\theta_{ij}$ determine homotopies of the loops.
We are thus in the position of having a map $\partial\Delta^2_+\wedge\Delta^1_+\rightarrow BQE$. Use the horn $\Lambda^1_0\subseteq\Delta^1$ to create a lift $\Delta^2_+\wedge\Delta^1_+\rightarrow BQE$.
Then, $\Delta^2_+\wedge\partial^0\Delta^1_+\rightarrow BQE$ is an element of $\pi_2 BQE$, say $\sigma$.

\begin{lemma}\label{lem:pi2}
    The element $\sigma\in\pi_2 BQE$ is the same as the class of $\pi_2 BQE$ canonically associated to the automorphism $\theta_{20}\circ\theta_{12}\circ\theta_{01}$ of $M_0$.
    \begin{proof}
        Each face $\partial^n\Delta^2_+\wedge\Delta^1$, i.e. each homotopy $\theta_{ij}$, corresponds to map
        \begin{equation*}
            \theta_{ij}:S^2-\left(D^1\vee D^1\right)\rightarrow BQE,
        \end{equation*}
        where the restriction of $\theta_{ij}$ to the first (resp. second) boundary disc is $M_i$ (resp. $M_j$).
        Let $X=S^2-\left(D^1\vee D^1\right)$. Then, together, $\theta_{01}$ and $\theta_{12}$ induce a map
        from the connected sum of $X$ with itself glued along $M_1$. This connected sum is itself homotopy equivalent to $X$. Throwing $\theta_{20}$ into the picture,
        we get two maps $X\rightarrow BQE$ which on one disc boundary are $M_0$ and on the second are $M_2$. Thus, they induce a map from the connected sum of $X$ with itself along the figure eight $S^1\vee S^1$.
        But, this connected sum is homotopy equivalent to $S^2$. The homotopy class of this map is $\sigma$.

        Now we check that $\sigma$ agrees with the automorphism of $M_0$ above.

        Recall that the category $QE$ consists of the same objects as $E$ but has as morphisms from $L$ to $N$ the collection of diagrams $L\lfib M\rcof N$, where $\lfib$
        and $\rcof$ denote admissible surjections and injections in $E$, modulo isomorphisms of such diagrams which are equalities on $L$ and $N$.

        An element $M$ of $E$ gives rise to an element $r_M$ in $\K_0(E)$.
        As a based loop in $BQE$, this is constructed as the composition of $i_M:0\twoheadleftarrow 0\hookrightarrow M$ with the inverse
        of $q_M:0\lfib M\rcof M$. See~\cite{srinivas_kbook}.

        As explained in~\cite[pp. 43-45]{srinivas_kbook}, an exact sequence $0\rightarrow L\rcof M\rfib N\rightarrow 0$ in $E$ corresponds to a \emph{choice} of 
        homotopy between $r_M$ and $r_L\cdot r_N$. This homotopy is constructed explicitly in \emph{ibid.} as a map from the $2$-sphere with a wedge sum of three discs removed,
        where the boundaries of these discs are the correctly oriented $r_*$. 

        In particular, given an isomorphism $\theta:L\riso M$, we get a map
        \begin{equation*}
            a_\theta:S^2-\left(D^1\vee D^1\vee D^1\right)\rightarrow BQE
        \end{equation*}
        with three discs removed. Since $\theta$ is an isomorphism, one of those discs can be filled in in $BQE$, and we get a map $X\rightarrow BQE$. So, the map from the punctured $S^2$ is just
        a proof of the equality $[L]=[M]$ in $\K_0(E)$.

        When $\theta:L\riso L$, the map $a_\theta$ is a map from $X$ to $BQE$ whose restrictions to the two boundary discs is the same. Thus, $a_\theta$ induces a map $S^2\rightarrow BQE$.
        This is the element of $\pi_2 BQE$ associated to $\theta$.

        In our situation of the $M_i$ and the $\theta_{ij}$, we have maps $a_{\theta_{ij}}:X\rightarrow BQE$. It is easy to see, looking at the construction of these maps that
        $a_{\theta_{12}\circ\theta_{01}}:X\rightarrow BQE$ is homotopic to gluing the two maps $a_{\theta_{01}}$ and $a_{\theta_{12}}$ along their common boundary circle $M_1$.
        Similarly, gluing $a_{\theta_{20}}$ to $a_{\theta_{12}\circ\theta_{01}}$ along $M_2$ is homotopic to
        \begin{equation*}
            a=a_{\theta_{20}\circ\theta_{12}\circ\theta_{01}}.
        \end{equation*}
        As identified above, since $a$ agrees on its two boundary circles, $a$ induces a map $S^2\rightarrow BQE$ that is homotopic to $\sigma$. But, this map is also the map
        associated to the automorphism $\theta_{20}\circ\theta_{12}\circ\theta_{01}$ of $M_0$, so the lemma follows.
    \end{proof}
\end{lemma}

\begin{construction}\label{const:reconstruction}
    We describe how to reconstruct the class $\alpha\in\Hoh^2(U_{\et},\Gm)$ from the stack $\StProjA$.
    The rank $1$ objects of $\StProjA$ form a $\Gm$-gerbe, $\StPicA$.
    Recall from \cite[Section~IV.3.4]{giraud_cohomologie_1971} or \cite[Theorem~5.2.8]{brylinski_loop_1993}
    the following procedure.
    First, one takes a cover $\mathcal{U}$ of $U$ such that there is an object
    of $\mathcal{L}_i\in\StPicA_{U_i}$ for all $i$. Second, choose isomorphisms
    \begin{equation*}
        \sigma_i:\Aut(\mathcal{L}_i)\riso\Gm|_{U_i},
    \end{equation*}
    Third, pick isomorphisms
    \begin{equation*}
        \theta_{ij}^{\alpha}:\mathcal{L}_i\riso\mathcal{L}_j
    \end{equation*}
    on a suitable refinement $1$-hypercover $\mathcal{V}_A\rightarrow\mathcal{U}_I\rightarrow U$. The composition
    \begin{equation*}
        \theta_{ki}^\delta\circ\theta_{jk}^\beta\circ\theta_{ij}^\alpha
    \end{equation*}
    is an element of $\Aut(\mathcal{L}_i)(Z_{ijk}^{\alpha\beta\delta})$, where
    \begin{equation*}
        Z_{ijk}^{\alpha\beta\delta}=V_{ij}^{\alpha}\times_{U_j}V_{jk}^{\beta}\times_{U_k}V_{ik}^{\delta}.
    \end{equation*}
    Then, 
    \begin{equation*}
        \sigma_i(\theta_{ki}^{\delta}\circ\theta_{jk}^{\beta}\circ\theta_{ij}^\alpha)\in\Gm(Z_{ijk}^{\alpha\beta\delta})
    \end{equation*}
    defines a $2$-cocycle in $\Gm$ which is in the same cohomology class as $\alpha$.
\end{construction}

\begin{theorem}\label{thm:d}
    Suppose that $U$ is geometrically connected and quasi-separated. Let $\alpha\in\Hoh^2(U_{\et},\Gm)$. Then, $d_2^{\alpha}([1])=\alpha$, through the natural isomorphism
    $$\Hoh^2(U_{\et},\mathcal{K}_1^\alpha)\riso\Hoh^2(U_{\et},\mathcal{K}_1)\riso\Hoh^2(U_{\et},\Gm).$$
    \begin{proof}
        Let $\alpha\in\Hoh^2(U_{\et},\Gm)$ be defined by a \Cech cocycle $\check{\alpha}\in\check{\Hoh}^2(\mathcal{U}^{\bullet},\Gm)$ for a hypercover $\mathcal{V}\rightarrow\mathcal{U}\rightarrow U$.
        Using Theorem~\ref{thm:morphism}, in order to prove the theorem, it suffices to prove that $d_2([1])=\check{\alpha}$ in the \Cech approximation spectral sequence associated to $\mathcal{U}^{\bullet}$.
        (The existence of such a hypercover is where the quasi-separated hypothesis is used; see \cite[Theorem~V.7.4.1]{sga4.2}.)
        Let $$Z_{ijk}^{\alpha\beta\delta}=V_{ij}^{\alpha}\times_{U_j}V_{jk}^{\beta}\times_{U_k}V_{ik}^{\delta}.$$
        We may represent $[1]\in\check{\Hoh}^0(\mathcal{U}^{\bullet},\mathcal{K}_0^\alpha)$ by an $\alpha$-twisted line bundle $\mathcal{L}_i$ on each $U_i$ of $\mathcal{U}^{\bullet}$.
        A homotopy from $\mathcal{L}_i$ to $\mathcal{L}_j$ is just an isomorphism $$\theta_{ij}^{\alpha}:\mathcal{L}_i|_{V_{ij}^{\alpha}}\riso\mathcal{L}_j|_{V_{ij}^{\alpha}},$$
        where such an isomorphism exists, up to possibly refining the hypercover $\mathcal{U}^{\bullet}$. Then, by Theorem~\ref{thm:cechd2} and Lemma~\ref{lem:pi2},
        the class of the automorphism
        $$\theta_{ki}^{\delta}\circ\theta_{jk}^{\beta}\circ\theta_{ij}^{\alpha}$$
        of $\mathcal{L}_i|_{Z_{ijk}^{\alpha\beta\delta}}$ in $\K_1^{\alpha}(Z_{ijk}^{\alpha\beta\delta})$ is the $Z_{ijk}^{\alpha\beta\delta}$-component of $d_2^{\alpha}([1])$.

        The $\alpha$-twisted line bundles $\mathcal{L}_i$ on $U_i$ induce pointwise weak equivalences of $K$-theory presheaves $\phi_i:\K^{\alpha}|_{U_i}\riso\K|_{U_i}$
        by tensor product with $\mathcal{L}_i^{-1}$.
        As shown in \cite{antieau_index_2009}, these local morphisms patch to create natural isomorphisms of $K$-theory sheaves
        \begin{equation*}
            \phi_i:\mathcal{K}^{\alpha}_{k}\riso\mathcal{K}_{k},
        \end{equation*}
        and of $K$-cohomology groups, in particular of
        \begin{equation*}
            \Hoh^2(U_{\et},\mathcal{K}_1^{\alpha})\riso\Hoh^2(U_{\et},\mathcal{K}_1). 
        \end{equation*}
        
        Define $\sigma_i$ by fixing an isomorphism
        \begin{equation*}
            \sigma_i:\mathcal{L}_i\otimes\mathcal{L}_i^{-1}\riso\Gm|_{U_i},
        \end{equation*}
        possibly refining $\mathcal{U}^{\bullet}$. Then, the diagram
        \begin{equation*}
            \begin{CD}
                \Aut(\mathcal{L}_i)     @>\sigma_{i,\ast}>>    \Gm|_{U_i}\\
                @VVV                                    @VVV\\
                \K^{\alpha}_1(U_i)      @>\phi_i>>      \K_1(U_i)
            \end{CD}
        \end{equation*}
        is commutative, where $\sigma_{i,\ast}$ is the natural isomorphism induced by $\sigma_i$.

        The diagram and Construction~\ref{const:reconstruction} imply that $d_2^{\alpha}([1])$ maps to the image of $\check{\alpha}$ in the map $\check{\Hoh}^2(\mathcal{U}^{\bullet},\Gm)\rightarrow\check{\Hoh}^2(\mathcal{U}^{\bullet},\K_1)$.
    \end{proof}
\end{theorem}

\begin{remark}
    There is a more sophisticated version of this theorem which uses the $K$-theory ring spectrum $\K$. Then the $\K^{\alpha}$ are module spectra over this ring spectrum.
    In this situation, the descent spectral sequence for $\K^{\alpha}$ is a module over the descent spectral sequence for $\K$, and thus, if $x\in\Hoh^s(U_{\et},\mathcal{K}_t^{\alpha})$,
    we may write $x=1_{\alpha}\cup y$, with $y\in\Hoh^s(U_{\et},\mathcal{K}_t)$ and $1\in\Hoh^0(U_{\et},\mathcal{K}_0^{\alpha})$ the generator. Thus,
    \begin{equation*}
        d_2^{\alpha}(x)=d_2^{\alpha}(1_{\alpha})\cup x\pm 1_{\alpha}\cup d_2(y)=\alpha\cup x\pm 1_{\alpha}\cup d_2(y).
    \end{equation*}
    So, all of the differentials in the $\Eoh_2$-page are determined by the theorem and the differentials of the descent spectral sequence for $\K$.
\end{remark}

\begin{theorem}[\textbf{Divisibility}]\label{thm:main}
    Let $U$ be geometrically connected and quasi-separated. If $\alpha\in\Hoh^2(U_{\et},\Gm)$, then
    \begin{equation*}
        per(\alpha)|eti(\alpha).
    \end{equation*}
    \begin{proof}
        This follows immediately from Theorem~\ref{thm:d}. Indeed, since $d_2^{\alpha}([1])=\alpha$, it follows that $per(\alpha)$
        generates $\Eoh_3^{0,1}\subseteq\Hoh^0(U_{\et},\ZZ)$ in the Brown-Gersten spectral sequence. Therefore, the generator of $\Eoh_{\infty}^{0,1}$ belongs to the subgroup $per(\alpha)\cdot\ZZ$.
        In other words, $per(\alpha)$ divides $eti(\alpha)$.
    \end{proof}
\end{theorem}

\begin{remark}
    In~\cite{kahn_levine}, there is an Atiyah-Hirzebruch spectral sequence in \'etale motivic cohomology
    \begin{equation*}
        \Hoh^{p-q}(X_{\et},\ZZ(-q))\Rightarrow\K^{\et}_{-p-q}(\mathcal{A}),
    \end{equation*}
    where $\mathcal{A}$ represents a class $\alpha\in\Br(X)$. Kahn and Levine show that
    \begin{equation*}
        d_2([1])=\alpha
    \end{equation*}
    in $\Hoh^3(X_{\et},\ZZ(1))=\Br(X)$.

    Similarly, in~\cite{atiyah_twisted_2006}, in the Atiyah-Hirzebruch spectral sequence
    \begin{equation*}
        \Hoh^p(X,\ZZ(q/2))\Rightarrow KU^{p-q}(X)_{\alpha}
    \end{equation*}
    converging to twisted topological $K$-theory, it is shown that
    \begin{equation*}
        d_3([1])=\alpha
    \end{equation*}
    in $\Hoh^3(X,\ZZ)$.
\end{remark}

\section{Index and Spectral Index}\label{sec:indeti}

Let $k$ be a field of finite etale cohomological dimension. Then, the main result of \cite{antieau_index_2009} is that
\begin{equation}\label{eq:etiperper}
    eti(\alpha)|per(\alpha)^{[\frac{d}{2}]},
\end{equation}
for all $\alpha\in\Br(k)$ having period not divisible by a few small primes.
Define $e(k)$ to be the smallest integer such that $$eti(\alpha)|per(\alpha)^{e(k)}$$ for all $\alpha\in\Br(k)$ whose period is not divisible by the characteristic of $k$.
Such an integer exists by~\cite[Theorem~6.10]{antieau_index_2009}.
Moreover, Theorem~\ref{thm:main} says that $e(k)\geq 1$. Recall, in the notation of the introduction, the group
\begin{equation*}
    F^{\alpha}=\K_0^{\alpha,\et}(k)^{(0)}/\K_0^{\alpha}(k)^{(0)}.
\end{equation*}
Its order is $\frac{ind(\alpha)}{eti(\alpha)}$.

\begin{definition}
    Say that the field $k$ has property $B_n$ if $F^{\alpha}$ is of order at most $per(\alpha)^n$ for all $\alpha\in\Br(k)$.
\end{definition}

\begin{proposition}
    Let $k$ be a field of finite etale cohomological dimension $d$, and let $e(k)$ be the integer defined above. Then, the period-index conjecture holds for $k$
    if $k$ has property $B_{d-1-e(k)}$. If the period-index conjecture holds for $k$, then $k$ has property $B_{d-2}$.
    \begin{proof}
        Fix $\alpha\in\Br(k)$. Then, $eti(\alpha)|per(\alpha)^{e(k)}$. If the order of $F^{\alpha}$ is at most $per(\alpha)^{d-1-e(k)}$, then, $ind(\alpha)|eti(\alpha)^{d-1-e(k)}$, whence the first statement.
        To prove the second statement, it suffices to note that $\frac{ind(\alpha)}{eti(\alpha)}$ is at most $\frac{per(\alpha)^{d-1}}{per(\alpha)}$,
        by the period-index conjecture and Theorem~\ref{thm:main}.
    \end{proof}
\end{proposition}

\begin{question}
    If the period-index conjecture hold for $k$ does it imply that $k$ has property $B_{d-1-e(k)}$?
\end{question}

\begin{remark}
    The proposition is interesting because it reveals a connection between the period-index conjecture and \'etale descent for (twisted) K-theory.
\end{remark}

\begin{remark}
    We do not know the integer $e(k)$ for any fields besides $2$-dimensional fields where the period-index conjecture is known, in which case $e(k)=1$. Our guess is that it
    is as big as possible, namely at least $[\frac{d}{2}]$. Some evidence for this may be that the fields such as $k=\CC((t_1))\cdots((t_d))$ have $ind(\alpha)|per(\alpha)^{[\frac{d}{2}]}$, so
    it is natural to ask whether $eti(\alpha)=ind(\alpha)$ for these fields. This is a line of inquiry we are currently pursuing.
\end{remark}

\bibliographystyle{amsplain-pdflatex}
\bibliography{brauer,mypapers}

\noindent
Benjamin Antieau
[\texttt{\href{mailto:antieau@math.ucla.edu}{antieau@math.ucla.edu}}]\\
Department of Mathematics\\
University of California Los Angeles\\
Box 951555\\
Los Angeles, CA 90095-1555\\
USA

\end{document}